\documentclass[12pt,a4paper,twoside,reqno]{amsart}   

\usepackage{latexsym}
\usepackage{amsmath}
\usepackage{amsfonts}
\usepackage{amssymb}
\usepackage{lscape}
\usepackage{boxedminipage}
\usepackage[matrix,arrow]{xy}
\usepackage{longtable}

\numberwithin{equation}{section}

\theoremstyle{plain}
\newtheorem{thm}{Theorem}[section]
\newtheorem{cor}[thm]{Corollary}

\theoremstyle{definition}
\newtheorem{Def}[thm]{Definition}
\newtheorem{rem}[thm]{Remark}

\newcommand{\R}{\mathrm{I\!R}}
\newcommand{\N}{\mathrm{I\!N}}

\newcommand{\Z}{\mathchoice {\hbox{$\sf\textstyle Z\kern-0.4em
Z$}}{\hbox{$\sf\textstyle Z\kern-0.4em Z$}}{\hbox{$\sf\scriptstyle
Z\kern-0.3em Z$}}{\hbox{$\sf\scriptscriptstyle Z\kern-0.2em Z$}}}

\newcommand{\Q}{\mathchoice {\setbox0=\hbox{$\displaystyle\rm
Q$}\hbox{\raise0.15\ht0\hbox to0pt{\kern0.4\wd0\vrule
height0.8\ht0\hss}\box0}}{\setbox0=\hbox{$\textstyle\rm
Q$}\hbox{\raise0.15\ht0\hbox to0pt{\kern0.4\wd0\vrule
height0.8\ht0\hss}\box0}}{\setbox0=\hbox{$\scriptstyle\rm
Q$}\hbox{\raise0.15\ht0\hbox to0pt{\kern0.4\wd0\vrule
height0.7\ht0\hss}\box0}}{\setbox0=\hbox{$\scriptscriptstyle\rm
Q$}\hbox{\raise0.15\ht0\hbox to0pt{\kern0.4\wd0\vrule
height0.7\ht0\hss}\box0}}}

\newcommand{\OO}{\mathchoice {\setbox0=\hbox{$\displaystyle\rm
O$}\hbox{\hbox to0pt{\kern0.4\wd0\vrule
height0.9\ht0\hss}\box0}}{\setbox0=\hbox{$\textstyle\rm O$}\hbox{\hbox
to0pt{\kern0.4\wd0\vrule
height0.9\ht0\hss}\box0}}{\setbox0=\hbox{$\scriptstyle\rm O$}\hbox{\hbox
to0pt{\kern0.4\wd0\vrule
height0.9\ht0\hss}\box0}}{\setbox0=\hbox{$\scriptscriptstyle\rm
O$}\hbox{\hbox to0pt{\kern0.4\wd0\vrule height0.9\ht0\hss}\box0}}}

\newcommand{\SL}{\mathrm{SL}}

\newcommand{\eps}{\varepsilon}
\newcommand{\vi}{\varphi}
\newcommand{\vkap}{\varkappa}

\newcommand{\qmq}[1]{\quad\mbox{#1}\quad}
\newcommand{\Menge}[2]{\{\,#1\,|\,#2\,\}}

\newcommand{\IM}{\mathop{\mathrm{Im}}\nolimits}

\newcommand{\tr}{\mathop{\mathrm{tr}}\nolimits}

\newcommand{\ord}{\mathop{\mathrm{ord}}\nolimits}

\newcommand{\As}{\mathrm{As}}

\newcommand{\Mengegr}[2]{\{\,#1\,{\bigr |}\,#2\,\}}

\newcommand{\wt}{\widetilde}

\newcommand{\C}{\mathchoice {\setbox0=\hbox{$\displaystyle\rm
C$}\hbox{\hbox to0pt{\kern0.4\wd0\vrule
height0.95\ht0\hss}\box0}}{\setbox0=\hbox{$\textstyle\rm C$}\hbox{\hbox
to0pt{\kern0.4\wd0\vrule
height0.95\ht0\hss}\box0}}{\setbox0=\hbox{$\scriptstyle\rm C$}\hbox{\hbox
to0pt{\kern0.4\wd0\vrule
height0.95\ht0\hss}\box0}}{\setbox0=\hbox{$\scriptscriptstyle\rm
C$}\hbox{\hbox to0pt{\kern0.4\wd0\vrule height0.95\ht0\hss}\box0}}}

\newcommand{\unity}{{1\!\!\!\:\mathrm{l}}}
\newcommand{\Pot}{\mathsf{Pot}}

\newcommand{\Div}{\mathsf{Div}}
\newcommand{\Jac}{\mathrm{Jac}}

\begin{document}

\title[Spectral data for simply periodic solutions of sinh-Gordon]{Spectral data for simply periodic solutions \\ of the sinh-Gordon equation}

\author[S.~Klein]{Sebastian Klein}

\address{School of Mathematical Sciences, University College Cork, Western Gateway Building, Western Road, Cork City, Ireland}
\email{s.klein@math.uni-mannheim.de}

\date{October 25, 2016.}



\abstract
{This note summarizes results that were obtained by the author in his habilitation thesis concerning the development of a spectral theory for 
simply periodic, 2-dimensional, complex-valued solutions \,$u$\, of the sinh-Gordon equation. Spectral data
for such solutions are defined for periodic Cauchy data on a line (following \textsc{Hitchin} and \textsc{Bobenko}) 
and the space of spectral data is described by an asymptotic characterization. Using methods of asymptotic estimates,
the inverse problem for the spectral data of such Cauchy data is answered. Finally a Jacobi variety for the spectral curve is constructed, 
and this is used to study the asymptotic behavior of the spectral data corresponding to actual simply periodic solutions of the sinh-Gordon equation
on strips of positive height.}

\endabstract

\maketitle

\bigskip

\tableofcontents

\section{Introduction}
\label{Se:intro}

The present paper constitutes a summary of the results obtained by the author in his habilitation thesis \cite{Klein:2015} concerning the development of a spectral theory for 
simply periodic, 2-dimensional, complex-valued solutions \,$u$\, of the sinh-Gordon equation. As such, it describes the constructions involved, and the most important
results, along with the fundamental ideas for their proofs. However, the detailed proofs of the results (some of which involve relatively lengthy calculations, e.g.~to obtain asymptotic
estimates) are referenced from \cite{Klein:2015}.

The primary object of the investigation are 
periodic, complex-valued solutions \,$u: X \to \C$\, of the 2-dimensional (i.e.~\,$X\subset \C$\,) sinh-Gordon equation
$$ \Delta u + \sinh(u)=0 \;. $$
We call such solutions \emph{simply periodic} when we wish to emphasize the difference to \emph{doubly periodic} solutions (which have two linear independent
periods). One important reason why solutions of the sinh-Gordon equation are interesting (apart from their relation to soliton theory) is that 
they arise from constant mean curvature surfaces without umbilical points in the 3-dimensional real space forms. 

\cite{Klein:2015} and the present paper set out to develop a spectral theory for simply periodic solutions of the sinh-Gordon equation.
The term ``spectral theory'' here refers ordinarily to a scheme of studying solutions of a given differential operator by looking at the spectrum of an associated Lax operator.
The eigenvalue equation for this operator can also be interpreted as the zero-curvature equation for a certain connection. 
The scheme of spectral theory
was first developed for the Korteweg-de Vries equation (KdV equation). A very accessible account of the spectral theory for the 1-dimensional Schr\"odinger operator
(which is the Lax operator for the KdV equation), has been given by \textsc{P\"oschel} and \textsc{Trubowitz} in \cite{Poeschel-Trubowitz:1987}.
This book has been very inspirational for my study of the spectral theory for the sinh-Gordon equation, and several results in Sections~\ref{Se:summary-asymp}
and \ref{Se:summary-remon} on the sinh-Gordon equation are analogous to corresponding results in \cite{Poeschel-Trubowitz:1987} for the 1-dimensional Schr\"odinger operator.

For the sinh-Gordon equation, the concept of a spectral theory is used in a somewhat more general sense however, in that one still considers the zero-curvature equation
for a certain (matrix-valued) connection, but the zero-curvature condition can no longer be interpreted as a eigenvalue equation.
We will still take the freedom to use the term spectral theory also in this case, and to apply the adjective ``spectral'' to the objects related to this theory.

The idea of a spectral theory for \emph{doubly} periodic solutions of the sinh-Gordon equation has first been applied by \textsc{Hitchin} in \cite{Hitchin:1990},
yielding a classification of the minimal tori in \,$S^3$\,. His results have 
later been refined by \textsc{Bobenko} and adapted to constant mean curvature immersions in all the 3-dimensional space forms. 
We mention that \textsc{Heller} has applied Hitchin's construction of spectral data to compact (closed) immersed surfaces of genus \,$g\geq 2$\, in \,$S^3$\,;
he obtains the most interesting results for surfaces which are ``Lawson symmetric'', i.e.~which have the symmetry group of one of the Lawson surfaces;
he also obtained constant mean curvature deformations of such surfaces in \,$S^3$\,. See for example \cite{Heller:2013}, \cite{Heller:2014} and \cite{Heller/Schmitt:2015}.
The present work differs from these previous results in that now \emph{simply} periodic solutions are considered, rather than doubly periodic solutions.

One of the most salient differences between the spectral theory for doubly periodic solutions and for simply periodic solutions of the sinh-Gordon equation
is that in the former case the spectral curve (a complex curve that comprises part of the spectral data for the sinh-Gordon equation)
is of finite geometric genus and can be compactified, whereas in the latter case, it generally
has infinite geometric genus. For this reason the classical results on compact Riemann surfaces, which were very useful for the study of doubly periodic solutions,
are not applicable in the present setting. We need to replace these results with specific arguments for open Riemann surfaces of the type of the spectral curves.
To make such arguments feasible, the behavior of the spectral curve and of the associated data near its ``open ends'' needs to be described, and this is the reason
why the asymptotic estimates for the spectral data play a very big role in the present study. 

Unfortunately there are only very few results on open Riemann surfaces with prescribed asymptotics found in the literature.
One example would be the book \cite{Feldman/Knoerrer/Trubowitz:2003} by \textsc{Feldman/Kn\"orrer/Trubowitz}. 
However, the results in the later part of the book, which would be very useful to us, depend on very strict geometric hypotheses for the surface under consideration, 
see \cite[Section~5]{Feldman/Knoerrer/Trubowitz:2003}, that are not satisfied for our spectral curves. 
For this reason we develop some results analogous to classical results on compact Riemann surfaces for spectral curves in this work, as needed.

In Section~\ref{Se:summary-spectrum}, we will construct spectral data for simply periodic solutions \,$u$\, of the sinh-Gordon equation, or more generally for 
so-called \emph{potentials}, i.e.~Cauchy data \,$(u,u_y)$\, for the sinh-Gordon equation, where \,$u$\, and \,$u_y$\, are periodic functions that are defined
only on a horizontal line. The spectral data consist of a complex curve \,$\Sigma$\,, called the \emph{spectral curve}, and a positive divisor \,$D$\, on \,$\Sigma$\,,
called the \emph{spectral divisor}. While it is possible for the spectral curve to have singularities, we will neglect the complications caused by them in 
the present summary; refer to \cite{Klein:2015} for their treatment. 

For the construction, we are interested in requiring only as weak regularity conditions for \,$(u,u_y)$\, as possible. 
There are two reasons: First, we are interested in characterizing precisely which divisors on a spectral curve are spectral divisors of some Cauchy data \,$(u,u_y)$\,;
it turns out that every additional differentiability condition imposed on \,$(u,u_y)$\, reduces the space of divisors by an intricate relationship between its divisor points.
By not imposing more regularity than necessary, we obtain a description of the space of divisors that is as simple as possible. 
Second, while any solution \,$u$\, of the sinh-Gordon equation is infinitely differentiable (in fact even real analytic, because the sinh-Gordon equation is elliptic)
on the interior of its domain, we are also interested in the behavior of the solution on the boundary of its domain, where its behavior can be worse. 
For these reasons we only require \,$(u,u_y) \in W^{1,2}([0,1]) \times L^2([0,1])$\,. 

In Section~\ref{Se:summary-asymp} we will describe the asymptotic behavior of the spectral data for a potential \,$(u,u_y)$\,. This information is fundamental 
for all following results. 

The \emph{inverse problem} for spectral data is the question if the solution \,$u$\, from which the spectral data is derived, or some other quantity associated
to \,$u$\,, is determined uniquely by the spectral data, and how it can be reconstructed from the spectral data. In Section~\ref{Se:summary-remon} we solve
the inverse problem for the monodromy \,$M(\lambda)$\, (defined in Section~\ref{Se:summary-asymp}); it turns out that the holomorphic functions comprising
the monodromy can be reconstructed explicitly as infinite sums and products in terms of the spectral data. 
After we have shown in Section~\ref{Se:summary-finite}
that the finite type spectral data are dense in the space of all spectral data (satisfying the asymptotic properties from Section~\ref{Se:summary-asymp}), 
we are able to solve the inverse problem for actual potentials \,$(u,u_y)$\, in Section~\ref{Se:summary-inverse}.

It remains to investigate the spectral data for actual simply periodic solutions \,$u$\, on horizontal strips of positive height, and to do so, we study
the flow of the spectral data under translations orthogonal to the direction of the period. For this purpose, we construct a Jacobi variety and an Abel map
for the spectral curve in Section~\ref{Se:summary-jacobi}. Like in the well-understood case of solutions of finite type, it turns out in Section~\ref{Se:summary-xytrans} that the motion
of the spectral divisor under translations is linear in the Jacobi coordinates. Using this result, we are finally able to describe
the asymptotic behavior of the spectral data \,$(\Sigma,D)$\, for actual simply periodic solutions \,$u$\, of the sinh-Gordon equation.
It turns out that they satisfy an exponential asymptotic law, much steeper than the asymptotic behavior of the spectral
data for Cauchy data \,$(u,u_y)$\,, as is to be expected, solutions of the sinh-Gordon equation being real analytic in the interior
of their domain.

\bigskip

\paragraph{\textbf{Acknowledgements.}} I would like to express my sincerest gratitude to Professor Martin Schmidt, who has advised me during the
creation of the underlying thesis \cite{Klein:2015}. His steady support and help has been invaluable to me. I have learned a lot from him. 
I would also like to thank Prof.~C.~Hertling, Dr.~A.~Klauer and Dr.~M.~Knopf for helpful discussions and advice.

\section{Spectral data for simply periodic solutions of the sinh-Gordon equation}
\label{Se:summary-spectrum}

Suppose that \,$X$\, is a horizontal strip in the complex plane \,$\C$\, with \,$0\in X$\, and that \,$u: X\to \C$\, is a (real or complex) solution
of the 2-dimensional sinh-Gordon equation
\begin{equation*}
\Delta u + \sinh(u)=0 \;
\end{equation*}
which is simply periodic with the period \,$1$\, in the sense that we have
$$ u(z+1)=u(z) \qmq{for all \,$z\in X$\,.} $$

We associate to \,$u$\, the family of linear partial differential equations \,$\mathrm{d} F_\lambda = \alpha_\lambda\cdot F_\lambda$\,
parameterized by the \emph{spectral parameter} \,$\lambda\in \C^*$\,, where the connection 1-form \,$\alpha_\lambda$\, is given by
\begin{align}
\alpha_\lambda 
& : = \frac{1}{4} \begin{pmatrix} i\,u_y & -e^{u/2} - \lambda^{-1}\,e^{-u/2} \\ e^{u/2}+\lambda\,e^{-u/2} & -i\,u_y \end{pmatrix} \mathrm{d}x \notag \\
\label{eq:summary-spectrum:alphaxy}
& \qquad\qquad + \frac{i}{4} \begin{pmatrix} -u_x & e^{u/2}-\lambda^{-1}\,e^{-u/2} \\ e^{u/2}-\lambda\,e^{-u/2} & u_x \end{pmatrix}\mathrm{d}y \; . 
\end{align}
The integrability condition for each of these partial differential equations is the Maurer-Cartan equation
\,$\mathrm{d}\alpha_{ \lambda} + [\alpha_{ \lambda} \wedge \alpha_{ \lambda}]=0$\, for \,$\alpha_\lambda$\,, which turns out to be equivalent to the 
sinh-Gordon equation for \,$u$\,. Therefore the differential equation \,$\mathrm{d} F_\lambda = \alpha_\lambda\cdot F_\lambda$\, is for every \,$\lambda\in \C^*$\,
integrable along any period of \,$u$\,. We denote the corresponding monodromy with base point \,$z_0=0$\,,
i.e.~for the integration along the interval \,$[0,1]\subset X$\,,
by \,$M(\lambda):=F_\lambda(1) \cdot F_\lambda(0)^{-1}$\,. In this way we obtain the monodromy map \,$M: \C^* \to \SL(2,\C),\; \lambda\mapsto M(\lambda)$\,.
\,$M(\lambda)$\, depends holomorphically on \,$\lambda$\,. 

We use the monodromy map to construct spectral data for the simply periodic solution \,$u$\,. The holomorphic function \,$\Delta := \tr M(\lambda): \C^*\to \C$\,
characterizes the complex curve defined by the eigenvalues of the monodromy, which we call the \emph{spectral curve}:
\begin{align}
\Sigma & := \Mengegr{(\lambda,\mu) \in \C^* \times \C}{\det(M(\lambda)-\mu\cdot \unity)=0} \notag \\
\label{eq:spectral:Sigma}
& = \Mengegr{(\lambda,\mu) \in \C^* \times \C}{\mu^2 - \Delta(\lambda)\cdot \mu + 1 = 0} \; . 
\end{align}
Because \,$\det M(\lambda)=1$\, holds for all \,$\lambda\in\C^*$\,, \,$\Sigma$\, is hyperelliptic by virtue of the holomorphic involution
$$ \sigma: \Sigma \to \Sigma,\; (\lambda,\mu) \mapsto (\lambda,\mu^{-1}) \; . $$
The branch points of this hyperelliptic curve are the zeros of \,$\Delta^2-4$\, of odd order. Note that it is possible for \,$\Sigma$\, to have singularities,
they occur exactly at the zeros of \,$\Delta^2-4$\, of order \,$\geq 2$\,. 

The spectral curve does not fully determine the monodromy \,$M(\lambda)$\, because it describes only its eigenvalues, not the corresponding eigenvectors.
The bundle of eigenvectors \,$\Lambda$\, 
of \,$M(\lambda)$\, on \,$\Sigma$\, is a holomorphic line bundle at least on \,$\Sigma'$\,, the Riemann surface of regular points of \,$\Sigma$\,. 
In general, such a line bundle is described by a divisor on \,$\Sigma$\,, but if \,$\Sigma$\, has singular points, then the concept of a divisor
is not so clear. The proper concept of divisor to use in this case is that of a generalized divisor introduced by \textsc{Hartshorne} in \cite{Hartshorne:1986},
i.e.~a subsheaf of the sheaf of meromorphic functions on \,$\Sigma$\, that is finitely generated over the sheaf of holomorphic functions on \,$\Sigma$\,.
For a detailed investigation of the structure of the divisor of \,$\Lambda$\, in the generalized sense, see \cite[Section~3]{Klein:2015}.
Note however that by applying a meromorphic transformation on the spectral curve, it is always possible to move the points in the support of the divisor of \,$\Lambda$\,
to regular points of \,$\Sigma$\,. For this reason, we will take the point of view throughout most of this summary that we consider
only those solutions \,$u$\, for which the support of the divisor
of \,$\Lambda$\, contains only regular points of \,$\Sigma$\,. 

Under this hypothesis, the eigenvector bundle \,$\Lambda$\, is described by a divisor \,$D$\, in the classical sense,
which we call the \emph{spectral divisor}. If we write \,$M(\lambda)  = \left( \begin{smallmatrix} a(\lambda) & b(\lambda) \\ c(\lambda) & d(\lambda) \end{smallmatrix} \right)$\,
with the holomorphic functions \,$a,b,c,d: \C^* \to \C$\,, then \,$\Sigma\to\C^2,\; (\lambda,\mu)\mapsto ( \tfrac{\mu-d(\lambda)}{c(\lambda)},1)$\, 
is a global meromorphic section of \,$\Lambda$\,, whence it follows that the spectral divisor is the polar divisor of the meromorphic function \,$\tfrac{\mu-d}{c}$\, 
on \,$\Sigma$\,. One can show that if a regular point \,$(\lambda_*,\mu_*)\in\Sigma$\, is in the support of \,$D$\,, hence a pole of \,$\tfrac{\mu-d}{c}$\, of some order \,$m$\,, then 
\,$c$\, has a zero of order exactly \,$m$\, and \,$\mu-a$\, has a zero of order at least \,$m$\, at that point. In particular we have \,$\mu_*=a(\lambda_*)$\,. 
It follows that the support of \,$D$\, consists of exactly those points \,$(\lambda_*,\mu_*) \in \Sigma$\, with \,$c(\lambda_*)=0$\, and \,$\mu=a(\lambda_*)$\,;
the multiplicity of such a point in \,$D$\, is given by the order of the zero of \,$c$\, at \,$\lambda=\lambda_*$\,. 

The spectral curve \,$\Sigma$\, and the spectral divisor \,$D$\, comprise the \emph{spectral data} for the simply periodic solution \,$u$\, of the sinh-Gordon equation. 

\bigskip

\textbf{Example.}
Let us look at the \emph{vacuum}, i.e.~the most obvious simply periodic solution \,$u=0$\, of the sinh-Gordon equation. It corresponds to a minimal surface of
zero sectional curvature, i.e.~to a minimal cylinder. 
The spectral data for the vacuum are of importance
because we will describe the asymptotic behavior of the spectral data (for \,$\lambda\to\infty$\, and \,$\lambda\to 0$\,) for general \,$u$\, by comparing the 
general spectral data to the spectral data of the vacuum. 

For \,$u=0$\, we obtain from Equation~\eqref{eq:summary-spectrum:alphaxy}
\begin{equation*}
\alpha_0 = \frac{1}{4} \begin{pmatrix} 0 & -(1+\lambda^{-1}) \\ 1+\lambda & 0 \end{pmatrix}\mathrm{d}x + \frac{i}{4} \begin{pmatrix} 0 & 1-\lambda^{-1} \\ 1-\lambda & 0 \end{pmatrix} \mathrm{d}y \; .
\end{equation*}
Because \,$\alpha_0$\, thus does not depend on \,$x$\,, we can calculate the monodromy of the vacuum simply as
\,$M_0(\lambda) = \exp \left( \tfrac{1}{4} \left( \begin{smallmatrix} 0 & -(1+\lambda^{-1}) \\ 1+\lambda & 0 \end{smallmatrix} \right) \right)$\,, and from there we obtain
\begin{equation*}
M_0(\lambda) = \begin{pmatrix} \cos(\zeta(\lambda)) & -{\lambda}^{-1/2}\,\sin(\zeta(\lambda)) \\ {\lambda}^{1/2}\,\sin(\zeta(\lambda)) & \cos(\zeta(\lambda)) \end{pmatrix} 
=: \begin{pmatrix} a_0(\lambda) & b_0(\lambda) \\ c_0(\lambda) & d_0(\lambda) \end{pmatrix}
\end{equation*}
with
\begin{equation}
\label{eq:summary-vacuum:zeta}
\zeta(\lambda) := \frac14\,\left({\lambda}^{1/2} + {\lambda}^{-1/2}\right) \; .
\end{equation}
Note that all the entries of \,$M_0$\, are even in \,${\lambda}^{1/2}$\,, and therefore indeed define holomorphic functions in \,$\lambda\in\C^*$\,.
We will use the names \,$a_0,\dotsc,d_0$\, for the component functions of the monodromy of the vacuum throughout the entire paper without any further reference,
and likewise
\begin{equation*}
\Delta_0(\lambda) := \tr(M_0(\lambda)) = 2\cos(\zeta(\lambda)) \; .
\end{equation*}

It follows that the spectral curve \,$\Sigma_0$\, of the vacuum is given by
\begin{align}
\Sigma_0 & = \Mengegr{(\lambda,\mu)\in \C^* \times \C}{\mu = \tfrac12\left( \Delta_0(\lambda) \pm \sqrt{\Delta_0(\lambda)^2-4} \right)} \notag \\
\notag
& = \Mengegr{(\lambda,\mu)\in \C^* \times \C}{\mu = e^{\pm i\,\zeta(\lambda)} } \; . 
\end{align}
This curve has no branch points above \,$\C^*$\,. It has double points at all those \,$\lambda\in \C^*$\, for which
\,$\zeta(\lambda)$\, is an integer multiple of \,$\pi$\,; these values of \,$\lambda$\, are exactly the following:
\begin{equation}
\label{eq:vacuum:lambdak0-def}
\lambda_{k,0} := 8\pi^2k^2 + 4\pi k \sqrt{4\pi^2k^2-1} - 1 \qmq{with \,$k\in \Z$\,.}
\end{equation}

We have \,$\lambda_{k,0}\in \R$\, for all \,$k\in \Z$\,, moreover we have the following asymptotic estimates for \,$\lambda_{k,0}$\,:
\begin{align}
\lambda_{k,0} & = 16\pi^2 k^2 - 2 + O(k^{-2}) \text{ for \,$k\to\infty$\,} \notag \\
\notag
\qmq{and} \lambda_{k,0} & = \frac{1}{16\pi^2}k^{-2} + \frac{1}{128\pi^4}k^{-4} + O(k^{-6}) \text{ for \,$k\to -\infty$\,.} 
\end{align}
In particular, \,$\lambda_{k,0}$\, tends to \,$\infty$\, resp.~to \,$0$\, for \,$k\to\infty$\, resp.~\,$k\to-\infty$\,.

The points \,$(\lambda_*,\mu_*)$\, in the support of the spectral divisor \,$D_0$\, of the vacuum are exactly those points for which \,$\lambda_*$\, is a zero of \,$c_0$\,
and \,$\mu=a_0(\lambda_*)$\, holds. It turns out that \,$c_0(\lambda_*)=0$\, holds if and only if \,$\lambda_*=\lambda_{k,0}$\, holds for \,$k\in \Z$\,; all these zeros
are of order \,$1$\,. Moreover we have \,$a_0(\lambda_{k,0})=(-1)^k =: \mu_{k,0}$\,. Therefore the divisor \,$D_0$\, of the vacuum is given by its support
$$ \Menge{(\lambda_{k,0},\mu_{k,0})}{k\in \Z} \;; $$
all these points have multiplicity \,$1$\, in \,$D_0$\,. 

Note that the spectral data of the vacuum do not satisfy our general hypothesis that the divisor points are regular points of the spectral curve; rather all divisor points
of the vacuum lie in double points of the corresponding spectral curve. 

\section{Asymptotic behavior of the spectral data}
\label{Se:summary-asymp}

In the present and in the following section, we suppose that only \emph{Cauchy data} for the periodic solution \,$u$\, are given. This means that we suppose that
we are given two functions \,$u$\, and \,$u_y$\, defined only on one period, namely on the real interval \,$[0,1]$\,. 
As explained in the Introduction, we want the differentiability condition for the Cauchy data to be as weak as possible. More specifically we require
that \,$u$\, is in the Sobolev space \,$W^{1,2}([0,1])$\, of weakly once-differentiable functions with square-integrable derivative. 
We require \,$u_y$\, only to be square-integrable, i.e.~that \,$u_y \in L^2([0,1])$\, holds.
We define ``mixed derivatives'' of \,$u$\, by using both \,$u$\, and \,$u_y$\,,
e.g.~we define \,$u_z = \tfrac12(u_x-iu_y)$\, where \,$u_x$\, is the Sobolev derivative of \,$u$\, and \,$u_y$\, is the function from the Cauchy data.

We are interested in Cauchy data that are periodic. Because of \,$u\in W^{1,2}([0,1])$\,, \,$u$\, is in particular continuous, 
so individual function values \,$u(x)$\, for \,$x \in [0,1]$\, are well-defined. The periodicity condition for the function \,$u$\,, which is defined at first 
only on \,$[0,1]$\, is then simply \,$u(0)=u(1)$\,. Note that there is no similar condition for \,$u_y$\,, because \,$u_y$\, is only square-integrable,
and therefore defined only up to null sets. 
We then regard \,$u$\, and \,$u_y$\, as being extended periodically to the real line. 

In the sequel, we will call such pairs \,$(u,u_y)$\, \emph{(periodic) potentials}. 
We denote the space of these potentials by
$$ \Pot := \Menge{(u,u_y) \in W^{1,2}([0,1]) \times L^2([0,1])}{u(0)=u(1)} \; . $$
\,$\Pot$\, becomes a Hilbert space via the inner product
$$ \langle (u,u_y)\;,\; (\wt{u},\wt{u}_y)\rangle_{\Pot} := \langle u,\wt{u} \rangle_{W^{1,2}([0,1])} + \langle u_y,\wt{u}_y \rangle_{L^2([0,1])} 
\qmq{for \,$(u,u_y), (\wt{u},\wt{u}_y) \in \Pot$\,.} $$

We can write down the \,$\mathrm{d}x$-part of the connection 1-form \,$\alpha_\lambda$\,, see Equation~\eqref{eq:summary-spectrum:alphaxy}, also for 
such potentials \,$(u,u_y) \in \Pot$\,:
$$ \alpha_\lambda = \frac{1}{4} \begin{pmatrix} i\,u_y & -e^{u/2} - \lambda^{-1}\,e^{-u/2} \\ e^{u/2}+\lambda\,e^{-u/2} & -i\,u_y \end{pmatrix} \mathrm{d}x \;, $$
and therefore the construction of spectral data for actual simply periodic solutions \,$u$\, of the sinh-Gordon equation from Section~\ref{Se:summary-spectrum}
carries over to periodic potentials. Thus we obtain a spectral family of monodromies \,$M(\lambda): \C^* \to \SL(2,\C)$\, and thereby spectral data
\,$(\Sigma,D)$\, for any given potential \,$(u,u_y) \in \Pot$\,.

As was explained in the Introduction, the asymptotic behavior of the monodromy and of the spectral data
for \,$\lambda\to\infty$\, and for \,$\lambda\to 0$\, is one of the fundamental tools in the present approach to simply periodic solutions of the sinh-Gordon equation.
To describe this asymptotic behavior for the monodromy \,$M(\lambda)$\,, we introduce certain spaces of holomorphic functions on \,$\C^*$\, whose members are
characterized by a certain asymptotic descent rate towards zero for \,$\lambda\to\infty$\, and/or for \,$\lambda \to0$\,.

For this purpose we will consider the Hilbert space \,$\ell^2$\, of square-summable sequences \,$(a_k)_{k\in \Z}$\, indexed over the integers; of course,
we equip \,$\ell^2$\, with the Hilbert space norm \,$\|a_k\|_{\ell^2} := \left( \sum_{k\in \Z} |a_k|^2 \right)^{1/2}$\,. Moreover, for \,$n,m\in \Z$\,
we consider the Hilbert space \,$\ell^2_{n,m}$\, of sequences \,$(a_k)_{k\in \Z}$\, defined by the Hilbert norm
$$ \|a_k\|_{\ell^2_{n,m}} := \left( \sum_{k=-\infty}^{-1} |k^m\cdot a_k|^2 + |a_0|^2 + \sum_{k=1}^\infty |k^n\cdot a_k|^2 \right)^{1/2} \;. $$

We also define for \,$k\in \Z$\,, where \,$\zeta(\lambda)$\, is as in Equation~\eqref{eq:summary-vacuum:zeta}
\begin{equation*}
S_k := \begin{cases}
\Mengegr{\lambda\in\C^*}{(k-\tfrac12)\pi \leq |\zeta(\lambda)| \leq (k+\tfrac12)\pi, |\lambda| > 1} & \text{for \,$k>0$\,} \\
\Mengegr{\lambda\in\C^*}{|\zeta(\lambda)|\leq \tfrac{\pi}{2}} & \text{for \,$k=0$\,} \\
\Mengegr{\lambda\in\C^*}{(-k-\tfrac12)\pi \leq |\zeta(\lambda)| \leq (-k+\tfrac12)\pi, |\lambda| < 1} & \text{for \,$k<0$\,.}
\end{cases}
\end{equation*}
Note that each \,$S_k$\, is a topological annulus, the \,$S_k$\, cover all of \,$\C^*$\,, and that \,$\lambda_{k,0} \in S_k$\, holds for every  \,$k\in \Z$\,. 

We then say that a holomorphic function \,$f: \C^* \to \C$\, has \emph{\,$\ell^2_{n,m}$-asymptotic of type \,$s$\,} (where \,$n,m\in \Z$\, and \,$s \geq 0$\,)
if there exists a sequence \,$(a_k)_{k\in \Z} \in \ell^2_{n,m}$\, of non-negative numbers so that
\begin{equation}
\label{eq:summary-asymp:As}
\forall k \in \Z \; \forall \lambda \in S_k \; : \; |f(\lambda)| \leq a_k \cdot e^{s\cdot |\IM(\zeta(\lambda))|}
\end{equation}
holds. We call any such sequence \,$(a_k)$\, a \emph{bounding sequence} for \,$f$\,, and denote the Banach space of all \,$\ell^2_{n,m}$-asymptotic functions
by \,$\As(\C^*,\ell^p_{n,m},s)$\,. If the condition \eqref{eq:summary-asymp:As} holds only for \,$k\geq 0$\, resp.~only for \,$k\leq 0$\, (instead of for all
\,$k\in \Z$\,), we say that \,$f$\, is \emph{\,$\ell^2_n$-asymptotic of type \,$s$\, for \,$\lambda\to\infty$\,} resp.~\emph{\,$\ell^2_m$-asymptotic of type \,$s$\,
for \,$\lambda\to 0$\,}, and we denote the space of such functions by \,$\As_\infty(\C^*,\ell^p_n,s)$\, resp.~by \,$\As_0(\C^*,\ell^p_m,s)$\,. 

The following theorem, which is of fundamental importance for the entire work, compares the monodromy \,$M(\lambda)$\, of a given periodic potential to the monodromy
\,$M_0(\lambda)$\, of the vacuum as described in the Example of Section~\ref{Se:summary-spectrum}. 

\begin{thm}
\label{T:summary-asymp:M}
Let \,$(u,u_y)\in \Pot$\, be given and \,$M(\lambda) = \left( \begin{smallmatrix} a(\lambda) & b(\lambda) \\ c(\lambda) & d(\lambda) \end{smallmatrix} \right)$\,
be the monodromy associated to \,$(u,u_y)$\,. We put \,$\tau := e^{-u(0)/2}$\,. Then we have
\begin{enumerate}
\item \,$a-a_0 \in \As(\C^*,\ell^2_{0,0},1)$\,
\item \,$b-\tau^{-1}\,b_0 \in \As_\infty(\C^*,\ell^2_1,1)$\, \quad and \quad \,$b-\tau\,b_0 \in \As_0(\C^*,\ell^2_{-1},1)$\,
\item \,$c-\tau\,c_0 \in \As_\infty(\C^*,\ell^2_{-1},1)$\, \quad and \quad \,$c-\tau^{-1}\,c_0 \in \As_0(\C^*,\ell^2_{1},1)$\,
\item \,$d-d_0 \in \As(\C^*,\ell^2_{0,0},1)$\,.
\end{enumerate}
\end{thm}

This theorem is proved in \cite{Klein:2015} in several stages: First, a weaker ``basic'' version is shown in \cite[Section~5]{Klein:2015}, 
where in the place of the \,$\ell^2$-sequences in 
Theorem~\ref{T:summary-asymp:M} one only has sequences which converge to zero. The proof of this basic version is based on a certain regauging of \,$\alpha$\, which
makes the leading term of \,$\alpha$\, (with respect to \,$\lambda$\,) independent of \,$u$\,, and thus makes an asymptotic estimate feasible. The proof of the basic version continues by
expressing the regauged monodromy as a power series in \,$u$\,, estimating the higher order terms of this power series, and eventually applying Riemann-Lebesgue's theorem.
A different version of the asymptotic estimate for \,$M(\lambda)$\, is shown in \cite[Section~7]{Klein:2015}: \,$M(\lambda_{k,0})$\, is asymptotically close to the Fourier coefficients
of \,$u_z$\, resp.~\,$-u_{\overline{z}}$\, (multiplied with certain powers of \,$k$\,). Because the Fourier coefficients of these \,$L^2$-functions are square-summable,
we obtain \,$\ell^2$-estimates for \,$M-M_0$\,, but only at the special points \,$\lambda_{k,0}$\,. Nonetheless, by combining both these versions of asymptotic estimates
for \,$M(\lambda)$\, it is possible to obtain the asymptotic estimate for the spectral divisor given in Theorem~\ref{T:summary-asymp:spectral}(2) below 
(\cite[Section 6 and 8]{Klein:2015}),  and then the infinite
sum resp.~product formulae for the entries of the monodromy described in the following section of this paper (\cite[Section 10]{Klein:2015}). By the examination of these formulas
one is finally liberated from the special role of the points \,$\lambda_{k,0}$\, in the \,$\ell^2$-version of the asymptotics and thereby one obtains the final form
of the asymptotics given in Theorem~\ref{T:summary-asymp:M} (see \cite[Section~11]{Klein:2015}). 

\bigskip

From Theorem~\ref{T:summary-asymp:M} 
it follows that also the spectral data
for any given potential \,$(u,u_y)$\, are asymptotically close to the spectral data for the vacuum. This is detailed in the following theorem. 

We say that there exist \emph{asymptotically and totally} exactly \,$m$\, points in every \,$S_k$\, with a certain property (where \,$m\in \N$\,), 
if there exists some \,$N\in \N$\, so that \,$S_k$\, contains exactly \,$m$\, points with this property for every \,$k\in \Z$\, with \,$|k|> N$\,,
and moreover \,$\bigcup_{|k|\leq N} S_k$\, contains exactly \,$m\cdot (2N+1)$\, points with the property.

\begin{thm}
\label{T:summary-asymp:spectral}
Let \,$(u,u_y)\in \Pot$\, be given, \,$M(\lambda)$\, be the monodromy associated to \,$(u,u_y)$\,, \,$\Delta := \tr M(\lambda)$\,, and \,$(\Sigma,D)$\, be the spectral
data for \,$(u,u_y)$\,.
\begin{enumerate}
\item
The function \,$\Delta^2-4$\, has asymptotically and totally exactly two zeros in every \,$S_k$\, (counted with multiplicities). They are the branch points resp.~the
singularities of \,$\Sigma$\,. It is thus possible to enumerate the zeros of \,$\Delta^2-4$\, by two sequences \,$(\vkap_{k,1})_{k\in \Z}$\, and \,$(\vkap_{k,2})_{k\in \Z}$\,
such that \,$\vkap_{k,\nu} \in S_k$\, holds for \,$|k|$\, large and \,$\nu\in\{1,2\}$\,, and then we have
$$ \vkap_{k,\nu}-\lambda_{k,0} \;\in\; \ell^2_{-1,3} $$
for \,$\nu\in \{1,2\}$\,.
\item
There is asymptotically and totally exactly one point \,$(\lambda_*,\mu_*)$\, in the support of \,$D$\, with \,$\lambda_* \in S_k$\, (counted with multiplicity). It is thus possible 
to enumerate the support of \,$D$\, by a sequence \,$(\lambda_k,\mu_k)_{k\in \Z}$\, so that \,$\lambda_k \in S_k$\, holds for \,$|k|$\, large, and then we have
$$ \lambda_k - \lambda_{k,0} \in \ell^2_{-1,3} \qmq{and} \mu_k-\mu_{k,0} \in \ell^2 \; . $$
\end{enumerate}
\end{thm}

The proof of part (2) of the above theorem uses only the two preliminary versions of the asymptotics of the monodromy from \cite[Sections 5 and 7]{Klein:2015}, see 
\cite[Sections 6 and 8]{Klein:2015}. This fact is important because the asymptotics for the support of the divisor are used to derive the sum and product formulas 
which reconstruct \,$M(\lambda)$\, from the spectral divisor, which are in turn used, among other things, to obtain the full strength of the asymptotic estimate
for the monodromy in the form of Theorem~\ref{T:summary-asymp:M}. In contrast, the proof of part (1) of Theorem~\ref{T:summary-asymp:spectral} uses the full 
Theorem~\ref{T:summary-asymp:M}, see \cite[Section 11]{Klein:2015}.

\section{Reconstruction of the monodromy}
\label{Se:summary-remon}

The inverse problem consists in reconstructing the potential resp.~the solution of the sinh-Gordon equation from the spectral data \,$(\Sigma,D)$\,. 
The first step in solving the inverse problem is to obtain the monodromy function \,$M: \C^* \to \SL(2,\C)$\,  from the spectral data.
It turns out that this can be done in a fairly explicit way: The entries of the matrix-valued function \,$M$\,, seen as holomorphic functions \,$\C^*\to \C$\,,
can be expressed as infinite sums or products in terms of the coordinates of the points in the support of \,$D$\,. By the same approach one also obtains
a first small glimpse at \,$u$\, itself, namely the function value \,$u(0)$\, can be obtained by an explicit formula in terms of the divisor points.

It will turn out that as long as the spectral divisor \,$D$\, does not contain any points of multiplicity \,$\geq 2$\,, the spectral divisor alone
(regarded as a set of points in \,$\C^* \times \C^*$\,) already implicitly determines the spectral curve on which it lies uniquely. To facilitate formulating
this insight, we will regard (spectral) divisors \,$D$\, on a curve \,$\Sigma$\,  also as plain sets of points in \,$\C^* \times \C^*$\,  with multiplicities in 
the sequel.

\begin{Def}
\label{D:summary-remon:asymptotic-nonspecial}
Let \,$D$\, be a positive divisor, regarded as a set of points in \,$\C^* \times \C$\, with multiplicities. 
\begin{enumerate}
\item We say that \,$D$\, is \emph{asymptotic} if there exists a sequence \,$(\lambda_k,\mu_k)_{k\in \Z}$\, that enumerates the points of \,$D$\, with their
multiplicities, and so that 
\begin{equation}
\label{eq:summary-remon:asymp}
\lambda_k-\lambda_{k,0}\in \ell^2_{-1,3} \qmq{and} \mu_k-\mu_{k,0} \in \ell^2_{0,0}
\end{equation}
holds. 
\item If \,$D$\, is asymptotic, we say that \,$D$\, is \emph{non-special}, if for every \,$k,\wt{k} \in \Z$\, with \,$\lambda_k=\lambda_{\wt{k}}$\, we also
have \,$\mu_{k}=\mu_{\wt{k}}$\,. 
\end{enumerate}
\end{Def}

If \,$D$\, is the spectral divisor of a potential \,$(u,u_y)\in\Pot$\,, it follows from Theorem~\ref{T:summary-asymp:spectral}(2) that \,$D$\, is asymptotic, and
because the present summary operates under the general hypothesis that no points of \,$D$\, are in singularities of the underlying spectral curve, 
the support of \,$D$\, consists of those points \,$(\lambda_*,\mu_*)$\, for which \,$c(\lambda_*)=0$\, and \,$\mu_* = a(\lambda_*)$\, holds 
(see Section~\ref{Se:summary-spectrum}), whence we can conclude immediately that \,$D$\, is also non-special.

The following theorem describes the reconstruction of the monodromy \,$M(\lambda)$\, and also of the value of \,$u(0)$\, from the spectral data. 
In particular, the reconstruction of the function \,$c$\,, whose zeros are known by the \,$\lambda$-components of the points in the 
spectral divisor, is in a way an adaption of Hadamard's Factorization Theorem to the present situation. The most significant difference between Hadamard's Theorem
and our situation is that the former concerns entire functions with zeros accumulating at \,$\lambda=\infty$\,, whereas we are interested in holomorphic
functions on \,$\C^*$\, whose zeros accumulate near both \,$\lambda=\infty$\, and \,$\lambda=0$\,. Notice that akin to Hadamard's Factorization Theorem
we also obtain an explicit representation of \,$c$\, as an infinite product. 

\begin{thm}
\label{T:summary-remon:reconstruction}
Let \,$D$\, be a non-special, asymptotic divisor on a spectral curve \,$\Sigma$\,. 
Then \,$D$\, is enumerated by a sequence \,$(\lambda_k,\mu_k)_{k\in \Z}$\, with the asymptotic behavior of 
\eqref{eq:summary-remon:asymp}. We define \,$\tau \in \C^*$\, and holomorphic functions \,$a,b,c,d: \C^* \to \C$\, in the following way, where the involved
infinite products and sums converge absolutely:

\begin{enumerate}
\item[(i)]
$$ \tau = \left( \prod_{k\in \Z} \frac{\lambda_{k,0}}{\lambda_k} \right)^{1/2} \; . $$
\item[(ii)] 
$$ c(\lambda) = \frac14\,\tau\,(\lambda-\lambda_0) \cdot \prod_{k=1}^\infty \frac{\lambda_{k}-\lambda}{16\,\pi^2\,k^2} \cdot \prod_{k=1}^{\infty} \frac{\lambda-\lambda_{-k}}{\lambda} \; . $$
\item[(iii)]
If the \,$\lambda_k$\, are pairwise unequal (i.e.~\,$c$\, has no zeros of order \,$\geq 2$\,), then the function \,$a$\, is obtained by the simple formula
$$ a(\lambda) = \sum_{k\in \Z} \frac{\mu_k \cdot c(\lambda)}{c'(\lambda_k)\cdot (\lambda-\lambda_k)} \; . $$

If some of the \,$\lambda_k$\, are equal (i.e.~\,$c$\, has zeros of higher order), then we need a more complicated method to reconstruct \,$a$\,:
For all \,$k\in \Z$\, let \,$d_k := \#\Menge{\wt{k}\in \Z}{\lambda_{\wt{k}}=\lambda_k} = \ord_{\lambda_k}(c)$\, be the multiplicity of \,$\lambda_k$\, in the support of \,$D$\,,
and let \,$\Lambda := \Menge{k\in \Z}{d_k \geq 2}$\, be the set of the indices \,$k$\, of the \,$\lambda_k$\, of higher order. \,$\Lambda$\, is finite.
For each \,$k\in \Lambda$\,, \,$(\lambda_k,\mu_k)$\, cannot be a branch point of \,$\Sigma$\,,%
\footnote{This is true only under the general hypothesis that no divisor points occur in singularities of \,$\Sigma$\,, which we made at the beginning of this summary.
For the general case it is again necessary to regard the spectral divisor as a generalized divisor, for the details see \cite[Sections 3 and 12]{Klein:2015}.}
and therefore we 
 can regard \,$\mu$\, as a holomorphic function in \,$\lambda$\, on a neighborhood of \,$(\lambda_k,\mu_k)$\,. We then choose \,$t_{k,1},\dotsc,t_{k,d_k}\in \C$\, so that 
with 
$$ A_k(\lambda) := \sum_{j=1}^{d_k} t_{k,j}\cdot \frac{c(\lambda)}{(\lambda-\lambda_k)^j} $$
we have
$$ A_k^{(\ell)}(\lambda_k) = \mu^{(\ell)}(\lambda_k) \qmq{for \,$\ell \in \{0,\dotsc,d_k-1\}$\,.} $$
Then 
$$ a(\lambda) = \sum_{k\in \Z\setminus \Lambda} \frac{\mu_k \cdot c(\lambda)}{c'(\lambda_k)\cdot (\lambda-\lambda_k)} + \sum_{k\in \Lambda} \frac{1}{d_k}\,A_k(\lambda) \; . $$
\item[(iv)]
The function \,$d$\, is obtained analogously to \,$a$\,, with \,$\mu$\, and \,$\mu_k$\, replaced by \,$\mu^{-1}$\, and \,$\mu_k^{-1}$\,, respectively.
\item[(v)]
Finally \,$b$\, is determined in terms of the other functions by the equation \,$ad-bc=1$\,. 
\end{enumerate}
Then the holomorphic functions \,$a,b,c,d:\C^* \to \C$\, are uniquely characterized by the following two properties:
They have the asymptotic behavior described in Theorem~\ref{T:summary-asymp:M}(1)--(4), and
\,$D$\, is the divisor so that the \,$\lambda_k$\, are all the zeros of \,$c$\, (with multiplicity) and \,$\mu_k = a(\lambda_k)$\,. 

Moreover, if \,$D$\, is the spectral divisor of some potential \,$(u,u_y)\in \Pot$\,, then 
\,$M(\lambda):=\left( \begin{smallmatrix} a(\lambda) & b(\lambda) \\ c(\lambda) & d(\lambda) \end{smallmatrix} \right)$\, 
is the monodromy of \,$(u,u_y)$\,. Moreover \,$e^{-u(0)/2}=\tau$\, holds; the latter
formula uniquely determines \,$u(0)$\, up to an integer multiple of \,$2\pi i$\,. 
\end{thm}

The proof of this theorem is worked out in detail in \cite[Sections 6, 10 and 12]{Klein:2015}, but the ideas are as follows: 
One can prove that there exists up to sign at most one holomorphic function \,$c: \C^* \to \C$\, which satisfies the asymptotics of Theorem~\ref{T:summary-asymp:M}(3)
(with whatsoever value of \,$\tau$\,)
and which has zeros exactly at the \,$\lambda_k$\, (counted with multiplicity). Moreover, if such a function \,$c$\, exists, then a holomorphic function 
\,$a: \C^* \to \C$\, with the asymptotics of Theorem~\ref{T:summary-asymp:M}(1) is uniquely determined by prescribing for every \,$k \in \Z$\,
the values of \,$a(\lambda_k),a'(\lambda_k),\dotsc,a^{(d_k-1)}(\lambda_k)$\,. This shows that \,$a,\dotsc,d$\, are uniquely determined by the properties given in the theorem.

On the other hand, 
it follows from the asymptotic assessments \eqref{eq:summary-remon:asymp} that the infinite product defining \,$\tau$\, converges in \,$\C$\, absolutely, 
that \,$\Lambda$\, is finite
and that the infinite product resp.~sums defining \,$c$\,, \,$a$\, and \,$d$\, in the theorem converge absolutely and locally uniformly in \,$\lambda \in \C^*$\,, thus they
indeed define holomorphic functions \,$\C^* \to \C$\,. It is easy to check that the zeros of \,$c$\, are exactly the \,$\lambda_k$\, (counted with
multiplicity), and that the equation \,$a(\lambda_k)=\mu_k$\, holds of order at least \,$d_k$\,. Moreover one sees that \,$b=\tfrac{ad-1}{c}$\, is holomorphic on \,$\C^*$\,,
even at the zeros of \,$c$\,. 
By a more involved analysis of the asymptotic behavior
of the infinite sums and products involved in the definition of the holomorphic functions \,$a,b,c,d$\,, one obtains that these functions
satisfy the asymptotic properties of Theorem~\ref{T:summary-asymp:M}(1)--(4); this result is also ultimately derived from the asymptotic behavior of
\,$\lambda_k$\, and \,$\mu_k$\,. 

Under the hypothesis of this summary that no spectral divisor points occur in singularities of the spectral curve, the spectral divisor \,$D$\, of some \,$(u,u_y)\in \Pot$\, 
is related to the spectral monodromy \,$M(\lambda) = \left( \begin{smallmatrix} a(\lambda) & b(\lambda) \\ c(\lambda) & d(\lambda) \end{smallmatrix} \right)$\, of that
potential by the fact that for any point \,$(\lambda_*,\mu_*)$\, in the support of \,$D$\,, say of order \,$m$\,, the function \,$c$\, has a zero at \,$\lambda_*$\,
exactly of order \,$m$\,, and \,$\mu-a(\lambda)$\, has a zero at \,$(\lambda_*,\mu_*)$\, at least of order \,$m$\, 
(see Section~\ref{Se:summary-spectrum}). Because the functions \,$a,\dotsc,d$\, 
are uniquely determined by these properties and their asymptotic behavior (Theorem~\ref{T:summary-asymp:M}), it follows that the functions \,$a,\dotsc,d$\, defined in 
Theorem~\ref{T:summary-remon:reconstruction} are indeed the entries of the monodromy of \,$(u,u_y)$\,. Because \,$\tau$\, is uniquely determined by the
asymptotic behavior of \,$c$\,, we also obtain the formula \,$e^{-u(0)/2}=\tau$\,. 

\bigskip

Note that in Theorem~\ref{T:summary-remon:reconstruction}, the spectral curve \,$\Sigma$\, resp.~the holomorphic local function \,$\mu$\, is only used in the reconstruction
in the case that some of the \,$\lambda_k$\, are equal. Because the divisor \,$D$\, is non-special, this can happen only if \,$D$\, contains points with multiplicity
\,$\geq 2$\,. If this is not the case, then the functions \,$a,\dotsc,d: \C^* \to \C$\, are already uniquely determined by the point set \,$\mathrm{supp}(D)$\,. Because these
functions determine the spectral curve \,$\Sigma$\, by means of Equation~\eqref{eq:spectral:Sigma} via the function \,$\Delta = \tr M = a+d$\,, we arrive at the following
Corollary:

\begin{cor}
\label{C:summary-remon:spectral}
Let \,$(u,u_y)\in \Pot$\, be a potential and \,$(\Sigma,D)$\, be the corresponding spectral data. If \,$D$\, does not contain any points of multiplicity \,$\geq 2$\,,
then the set \,$\mathrm{supp}(D) \subset \C^* \times \C^*$\, already uniquely determines the spectral curve \,$\Sigma$\,. 
\end{cor}

\section{Divisors of finite type}
\label{Se:summary-finite}

Next we address the inverse problem for potentials, i.e.~we would like to show that for given data \,$(\Sigma,D)$\,, where \,$\Sigma$\,
is a hyperelliptic complex curve above \,$\C^*$\, with the asymptotic behavior of Theorem~\ref{T:summary-asymp:spectral}(1), and \,$D$\, is a non-special,
asymptotic divisor on \,$\Sigma$\,, there exists one and (essentially) only one potential \,$(u,u_y) \in \Pot$\, such that \,$(\Sigma,D)$\, are the spectral
data of \,$(u,u_y)$\,. 

To prove that this is indeed the case, we will use the fact that the inverse problem is already well-understood in a certain special case, namely the case
where the potential is of finite type. (Among the potentials of finite type are those which belong to doubly periodic solutions of the sinh-Gordon equation;
among them there are in turn the potentials corresponding to CMC tori, which have been classified by Pinkall/Sterling and Hitchin.) To be able to apply the already
known facts on finite type potentials, we show that the finite type potentials resp.~spectral data are dense in the space of all potentials resp.~spectral data.
One of course expects this result to be true from the experience with other integrable systems, but to my knowledge, no explicit proof for the case of the
sinh-Gordon integrable system is yet found in the literature. It turns out that a natural proof of this statement can be given in the context of the present paper.

\begin{Def}
Let \,$(\Sigma,D)$\, be given, where \,$\Sigma$\,
is a hyperelliptic complex curve above \,$\C^*$\, with the asymptotic behavior of Theorem~\ref{T:summary-asymp:spectral}(1), and \,$D$\, is a non-special,
asymptotic divisor on \,$\Sigma$\,.

We then say that \,$(\Sigma,D)$\, is of \emph{finite type} if the following two conditions hold:
\begin{enumerate}
\item
\,$\Sigma$\, has finite geometric genus (i.e.~only finitely many of the double points of the spectral curve of the vacuum have ``opened up'' into a pair of branch points
with positive distance).
\item
All but finitely many of the points in the support of \,$D$\, lie in double points of \,$\Sigma$\,. 
\end{enumerate}

We also say that a potential \,$(u,u_y)\in \Pot$\, is of \emph{finite type}, if the corresponding spectral data are of finite type. 
\end{Def}

Thus spectral data \,$(\Sigma,D)$\, look like the spectral data of the vacuum at all but finitely many of the divisor points. Note also that finite type
spectral data (like the vacuum spectral data) do not satisfy our general hypothesis that all divisor points are in regular points of the spectral data.
Strictly speaking, one would therefore need to consider generalized divisors to be able to handle finite type spectral data in our setting. 
For the purposes of the present summary paper, we will again ignore the technical complications arising from this fact, however. 

\begin{Def}
Let \,$D$\, be an asymptotic divisor and suppose that \,$\mathrm{supp}(D)$\, is enumerated by the sequence \,$(\lambda_k,\mu_k)_{k\in \Z}$\, in such a way that
\eqref{eq:summary-remon:asymp} holds. We say that \,$D$\, is \emph{tame} if the \,$\lambda_k$\, are pairwise unequal. 
\end{Def}

To simplify our construction, we restrict our consideration to tame divisors. Any tame divisor is necessarily non-special, and by Corollary~\ref{C:summary-remon:spectral} 
a tame divisor uniquely determines its spectral curve. When working with tame divisors it therefore suffices to consider the divisor itself as a point set in \,$\C^* \times \C^*$\,,
without considering the underlying spectral curve. It is clear that the set of tame divisors is open and dense in the space of all asymptotic divisors. To show that the finite type
spectral data are dense in the space of all spectral data, it therefore suffices to show that the finite type tame divisors are dense in the space of all tame asymptotic divisors.

\newpage

\begin{thm}
\label{T:summary-finite-type:dense}
The set of finite type spectral data is dense in the space of all spectral data.

More specifically, for any tame divisor \,$D$\,, enumerated by the sequence \,$(\lambda_k,\mu_k)_{k\in \Z}$\, with \eqref{eq:summary-remon:asymp}, and any \,$\eps>0$\,
there exists a tame divisor \,$D^*$\, of finite type, enumerated by the sequence \,$(\lambda_k^*,\mu_k^*)_{k\in \Z}$\, with \eqref{eq:summary-remon:asymp}, so that
\begin{equation}
\label{eq:summary-finite-type:estimate}
\|\lambda_k^*-\lambda_k\|_{\ell^2_{-1,3}} + \|\mu_k^*-\mu_k\|_{\ell^2_{0,0}} < \eps
\end{equation}
holds. Moreover for given \,$N\in \N$\,, the divisor \,$D^*$\,  can be chosen such that
$$ \forall \, k\in \Z, |k|\leq N \; : \; \lambda_k^* = \lambda_k,\; \mu_k^* = \mu_k $$
holds.
\end{thm}

The detailed proof of this theorem is found in \cite[Section 13]{Klein:2015}. It is based on an application of the Banach Fixed Point Theorem. 
We want to find a spectral curve \,$\Sigma^*$\, 
such that \,$(\lambda_k,\mu_k)\in \Sigma^*$\, holds for \,$|k|\leq N$\, and such that \,$\Sigma^*$\, has a double point
near \,$\lambda_{k,0}$\, for each \,$k$\, with \,$|k|>N$\,. The latter condition means: 
Denoting the trace function of \,$\Sigma^*$\, by \,$\Delta^*$\,,
there are zeros \,$(\eta_k)_{|k|>N}$\, of \,$\Delta'$\, for which
\,$\Delta(\eta_k)=2(-1)^k$\, holds. We seek to construct this trace function \,$\Delta^*$\, so that the corresponding spectral curve \,$\Sigma^*$\,
(defined by Equation~\eqref{eq:spectral:Sigma}) has the desired properties.

Because the \,$\lambda_k$\, are pairwise unequal (\,$D$\, being tame), one can show similarly as in the proof of the reconstruction of
the function \,$a$\, in Theorem~\ref{T:summary-remon:reconstruction} that for any sequence \,$(z_k)_{k\in \Z} \in \ell^2_{0,0}$\, there exists one and only one holomorphic function 
\,$\Delta: \C^*\to \C$\, with \,$\Delta-\Delta_0 \in \As(\C^*,\ell^2_{0,0},1)$\, and
$$ \Delta(\lambda_k) = 2(-1)^k + z_k \qmq{for all \,$k\in \Z$\,.} $$

We will use the Banach Fixed Point Theorem to determine \,$(z_k)$\, so that the corresponding function \,$\Delta$\, has the desired properties.
For this purpose we fix \,$N\in \N$\, and consider the Banach space \,$\mathfrak{B}_N$\, of 
square-summable sequences \,$(z_k)_{|k|>N}$\,, equipped with the \,$\ell^2$-norm. For given \,$(z_k) \in \mathfrak{B}_N$\, we construct another sequence \,$(\wt{z}_k) 
\in \mathfrak{B}_N$\, in the following way: Let \,$\Delta: \C^* \to \C$\, be the holomorphic function with \,$\Delta-\Delta_0 \in \As(\C^*,\ell^2_{0,0},1)$\, and
$$ \Delta(\lambda_k) = \begin{cases} \mu_k+\mu_k^{-1} & \text{for \,$|k|\leq N$\,} \\ 2(-1)^k + z_k & \text{for \,$|k|>N$\,} \end{cases} \; . $$
It can be shown that the zero set of \,$\Delta'$\, can be enumerated (with multiplicities) by a sequence \,$(\eta_k)_{k\in \Z}$\, so that \,$\eta_k-\lambda_{k,0}\in \ell^2_{-1,3}$\,
holds and then by one more zero \,$\eta_*\in \C^*$\,. We then define a new sequence \,$(\wt{z}_k)_{|k|>N}$\, by 
$$ \wt{z}_k := z_k-\bigr( \Delta(\eta_k)-2(-1)^k \bigr) \; . $$
It can be shown that \,$(\wt{z}_k) \in \ell^2_{0,0}$\, holds. Thus we define the iteration map \,$\Phi: \mathfrak{B}_N \to \mathfrak{B}_N,\; (z_k) \mapsto (\wt{z}_k)$\,. 

One then shows by a detailed analysis of the asymptotic behavior of all the quantities involved that \,$\Phi$\, is Lipschitz continuous on any closed ball in 
\,$\mathfrak{B}_N$\, and one also obtains an estimate for the value of the Lipschitz constant (in dependence on \,$N$\, and the radius of the ball). 
It follows from this investigation that if \,$N$\, is chosen large enough, and moreover some \,$\delta>0$\, is chosen small enough, 
then \,$\Phi$\, maps the closed ball of radius \,$\delta$\,
in \,$\mathfrak{B}_N$\, into itself, and is a contraction on that ball. It follows by the Banach Fixed Point Theorem that \,$\Phi$\, has one and only one fixed point
\,$(z_k^*)_{|k|>N}$\, in this ball. 

If we let \,$\Delta^*$\, and \,$(\eta_k^*)$\, be the objects defined above for this sequence \,$(z_k^*)$\,, then we have \,$\Delta^*(\eta_k^*) = 2(-1)^k$\, for \,$|k|>N$\,, and therefore
the spectral curve \,$\Sigma^*$\, corresponding to \,$\Delta^*$\, has a double point at \,$\eta_k^*$\, for every \,$k$\, with \,$|k|>N$\,. Therefore the divisor
\,$D^*$\, the support of which is given by \,$(\lambda_k^*,\mu_k^*)$\, with 
$$ \lambda_k^* := \begin{cases} \lambda_k & \text{for \,$|k|\leq N$\,} \\ \eta_k^* & \text{for \,$|k|>N$\,} \end{cases}
\qmq{and}
\mu_k^* := \begin{cases} \mu_k & \text{for \,$|k|\leq N$\,} \\ (-1)^k & \text{for \,$|k|>N$\,} \end{cases} $$
is asymptotic and of finite type, and satisfies \eqref{eq:summary-finite-type:estimate} provided that \,$N$\, is chosen large enough also in relation to \,$\eps$\,. 

\section{The inverse problem for potentials}
\label{Se:summary-inverse}

After we have shown that the finite type divisors are dense in all asymptotic divisors, we are now ready to 
 discuss the inverse problem for potentials. We  show that the potential \,$(u,u_y) \in \Pot$\, is uniquely determined by
its spectral divisor \,$D$\,, at least in the case where \,$D$\, is tame. 

To phrase this statement more precisely, let us denote by \,$\Div$\, the space of asymptotic divisors. In view of Definition~\ref{D:summary-remon:asymptotic-nonspecial}(1)
it seems tempting to identify \,$\Div$\, with the Banach space \,$\ell^2_{-1,3} \oplus \ell^2_{0,0}$\,. However, we need to be careful because the enumeration of the support
of \,$D$\, by a sequence \,$(\lambda_k,\mu_k)_{k\in \Z}$\, is only unique up to reordering finitely many of the elements. Thus we define a distance on \,$\Div$\, in
the following way: For given \,$D^{[1]},D^{[2]} \in \Div$\, with corresponding enumerations \,$(\lambda_k^{[\nu]},\mu_k^{[\nu]})_{k\in \Z}$\, of \,$\mathrm{supp}(D^{[\nu]})$\,
(\,$\nu\in\{1,2\}$\,) so that \eqref{eq:summary-remon:asymp} holds  we put
$$ \|D^{[1]}-D^{[2]}\|_{\Div} := \inf_{\sigma_1,\sigma_2 \in P(\Z)} \left( \left\|\lambda_{\sigma_1(k)}^{[1]}-\lambda_{\sigma_2(k)}^{[2]}\right\|_{\ell^2_{-1,3}}^2
+ \left\|\mu_{\sigma_1(k)}^{[1]}-\mu_{\sigma_2(k)}^{[2]} \right\|_{\ell^2_{0,0}}^2 \right)^{1/2} \; . $$
Here \,$P(\Z)$\, denotes the group of finite permutations of \,$\Z$\,, i.e.~of those permutations \,$\sigma: \Z\to\Z$\, for which the set \,$\Z \setminus \mathrm{Fix}(\sigma)$\,
is finite. 

We also consider the open and dense subset \,$\Div_{tame}$\, of tame divisors in \,$\Div$\,. Moreover we say that a potential \,$(u,u_y)$\, is \emph{tame} if the
corresponding spectral divisor is tame. We denote the subset of tame potentials in \,$\Pot$\, by \,$\Pot_{tame}$\,. It will turn out that also \,$\Pot_{tame}$\, is
an open and dense subset of \,$\Pot$\,. 

Our object of interest in the present section is the map \,$\Phi: \Pot \to \Div$\, that maps each potential \,$(u,u_y)$\, onto the corresponding spectral divisor \,$D$\,. 
The statement about the inverse problem is expressed by the following theorem:

\begin{thm}
\label{T:summary-inverse:diffeo}
\,$\Phi|\Pot_{tame}: \Pot_{tame}\to \Div_{tame}$\, is a diffeomorphism onto an open and dense subset of \,$\Div_{tame}$\,. 
\end{thm}

\begin{rem}
\,$\Phi$\, is not immersive at potentials \,$(u,u_y) \in \Pot\setminus \Pot_{tame}$\,. This is true even under our general hypothesis for this summary
that the spectral divisor corresponding
to \,$(u,u_y)$\, should not contain any singular points of \,$\Sigma$\,. Indeed, if a point occurs in the support of \,$D$\, with multiplicity \,$\geq 2$\,,
there is an 
entire family of integral curves of \,$x$-translation in \,$\Div$\, that intersect in \,$D$\,, and therefore \,$\Phi$\, cannot be immersive.
To make \,$\Phi$\, an immersion (and consequently a local diffeomorphism) near such points, we
would need to replace the range \,$\Div$\, of \,$\Phi$\, by a suitable blow-up at its singularities (see~e.g.~\cite{Hartshorne:1977}, p.~163ff.). 
We do not carry out such a construction here.

The reason why the image of \,$\Phi|\Pot_{tame}$\, is not all of \,$\Div_{tame}$\, is that even though any tame divisor \,$D\in\Div_{tame}$\, is non-special, it is possible for \,$D$\, to become special under 
\,$x$-translation. 
If this occurs, the potential corresponding to \,$D$\, has a singularity for the corresponding value of \,$x$\,, and thus \,$D$\, does not belong to a potential \,$(u,u_y)\in \Pot$\, in our sense.
The investigation of sinh-Gordon potentials with singularities, corresponding to divisors \,$D\in\Div$\, which become special under \,$x$-translation for some value of \,$x$\,, would be extremely interesting
in view of studying compact constant mean curvature surfaces.
\end{rem}

The proof of Theorem~\ref{T:summary-inverse:diffeo} is set out in \cite[Sections 14 and 15]{Klein:2015}. It is based on two different building blocks.

The first building block are the divisors of finite type. It has been shown by \textsc{Bobenko} \cite[Theorem~4.1]{Bobenko:1991a} 
(also compare the explicit construction in terms of vector-valued Baker-Akhiezer functions by \textsc{Knopf} in \cite[Proposition 4.34]{Knopf:2013})
that if \,$D$\, is an asymptotic divisor of finite type, 
so that the \,$x$-translation \,$D(x)$\, of \,$D$\, exists and is non-special for every \,$x\in [0,1]$\,, then there exists one and only one \,$(u,u_y)\in \Pot$\,
with \,$\Phi((u,u_y))=D$\,. 

Here we mean by the translation \,$D(x_0)$\, of the spectral divisor \,$D$\, corresponding to \,$(u,u_y)\in \Pot$\, 
by \,$x_0\in [0,1]$\, the spectral divisor of the translated potential \,$(u(x+x_0),u_y(x+x_0))\in \Pot$\,. The corresponding motion of the coordinates \,$\lambda_k,\mu_k$\, 
of the points in 
the support of \,$D$\, can be described by differential equations in the \,$\lambda_k$\, and \,$\mu_k$\,. Because a divisor point that is located in a double point
of the spectral curve does not move at all under translations, only finitely many coordinate functions are actually in motion under translation in the case
of a finite type divisor. Thus we can define  \,$D(x)$\, for divisors \,$D$\, of finite type at least for small \,$|x|$\, without reference
to a potential \,$(u,u_y)$\,. Note that if \,$D$\, is of finite type, then \,$D(x)$\, also is of finite type whenever \,$D(x)$\, is defined. 

Using the fact (Theorem~\ref{T:summary-finite-type:dense}) that the set of finite type tame divisors is dense in \,$\Div_{tame}$\,, one can also show that
the set \,$\Div_{ft,xtame}$\, of finite type divisors \,$D$\, so that \,$D(x)$\, is defined and tame for every \,$x\in [0,1]$\, is dense in \,$\Div_{tame}$\,. Because any 
tame divisor is non-special, the mentioned result by Bobenko implies that for any \,$D\in \Div_{ft,xtame}$\, there exists one and only one \,$(u,u_y)\in \Pot$\,
with \,$\Phi((u,u_y))=D$\,.

The second building block is a symplectic basis for the tangent space \,$T_{(u,u_y)}\Pot$\, for \,$(u,u_y) \in \Pot_{tame}$\,. The corresponding coordinates
are analogous to the coordinates on finite-dimensional symplectic spaces given by Darboux's Theorem, therefore we will call this basis \emph{Darboux coordinates} even
in the present, infinite-dimensional setting. We equip \,$T_{(u,u_y)}\Pot$\, with the non-degenerate symplectic form
\begin{equation}
\label{eq:darboux:Omega-defined}
\Omega: T_{(u,u_y)}\Pot \times T_{(u,u_y)}\Pot \to \C,\; \bigr( \, (\delta u,\delta u_y)\,,\, (\wt{\delta} u,\wt{\delta} u_y) \, \bigr) \mapsto \int_0^1 \bigr( \delta u \cdot \wt{\delta} u_y - \wt{\delta}u \cdot \delta u_y \bigr) \,\mathrm{d}x \; .
\end{equation}
Then it was shown by \textsc{M.~Knopf} in \cite{Knopf:2015}, together with the author in \cite[Section~14]{Klein:2015} that there exists a symplectic basis
\,$(v_k,w_k)_{k\in \Z}$\, of \,$T_{(u,u_y)}\Pot$\, with respect to \,$\Omega$\, which can be defined explicitly in terms of \,$u$\, and 
the ``extended frame'' of \,$(u,u_y)$\, (i.e.~of the solution \,$F_\lambda$\, of the partial differential equation \,$F_\lambda' = \alpha_\lambda\,F_\lambda$\,
with \,$F_\lambda(0)=\unity$\,). Moreover if we denote by \,$D$\, the spectral divisor of \,$(u,u_y)$\,, and enumerate its support by a sequence \,$(\lambda_k,\mu_k)_{k\in \Z}$\,
as in Definition~\ref{D:summary-remon:asymptotic-nonspecial}(1) once more, we regard \,$\lambda_k$\, and \,$\mu_k$\, as complex-valued functions defined at least
on a neighborhood of \,$(u,u_y)$\, in \,$\Pot_{tame}$\,. Then we can think of \,$\Phi$\, as the map \,$(u,u_y) \mapsto (\lambda_k,\mu_k)_{k\in \Z}$\,,
and the tangent space \,$T_D\Div$\, is spanned by the variations \,$\delta \lambda_k$\, and \,$\delta \mu_k$\,, where \,$k\in \Z$\,. In these terms
we define the non-degenerate symplectic form
$$ \wt{\Omega}: T_D\Div \times T_D\Div \to \C,\; \bigr(\, (\delta \lambda_k,\delta \mu_k)\,,(\wt{\delta} \lambda_k,\wt{\delta} \mu_k)\,\bigr) \mapsto 
\frac{i}{2}\,\sum_{k\in \Z} \left( \frac{\delta \lambda_k}{\lambda_k} \cdot \frac{\wt{\delta} \mu_k}{\mu_k} - \frac{\wt{\delta} \lambda_k}{\lambda_k} \cdot \frac{\delta \mu_k}{\mu_k} \right) $$
on \,$T_D\Div$\,. By Knopf and the author it was moreover shown at the cited locations that 
for \,$(\delta u,\delta u_y)\,,\,(\wt{\delta} u,\wt{\delta} u_y) \in T_{(u,u_y)}\Pot$\,, we have
\begin{equation}
\label{eq:summary-inverse:OmegaWtOmega}
\Omega\bigr( \, (\delta u,\delta u_y)\,,\,(\wt{\delta} u,\wt{\delta} u_y) \,\bigr) = \wt{\Omega}\bigr(\, (\delta \lambda_k,\delta \mu_k)\,,(\wt{\delta} \lambda_k,\wt{\delta} \mu_k)\,\bigr) \;,
\end{equation}
where \,$(\delta \lambda_k,\delta \mu_k)$\, is the variation of \,$(\lambda_k,\mu_k)$\, corresponding to \,$(\delta u,\delta u_y)$\,.

Using these two building blocks, one can prove Theorem~\ref{T:summary-inverse:diffeo}, i.e.~that \,$\Phi|\Pot_{tame}$\, is a diffeomorphism. We first note
that it is clear from the construction of the spectral data that \,$\Phi$\, is smooth in the ``weak'' sense that all the coordinate functions
\,$\lambda_k,\mu_k$\, of the spectral divisor are smooth near \,$(u,u_y)$\,. For given \,$(u,u_y) \in \Pot_{tame}$\, we thus have
a ``weak'' derivative of \,$\Phi$\, at \,$(u,u_y)$\,, namely the linear map \,$\Phi'((u,u_y)): T_{(u,u_y)}\Pot \to T_D \Div$\,, 
where \,$D := \Phi((u,u_y))$\,. It follows from Equation~\eqref{eq:summary-inverse:OmegaWtOmega} that \,$\Phi'((u,u_y))$\, is in fact a 
symplectomorphism between \,$(T_{(u,u_y)}\Pot,\Omega)$\, and \,$(T_D\Div,\wt{\Omega})$\,. 

To prove that \,$\Phi$\, is also differentiable at \,$(u,u_y)$\, in the stronger sense,
namely as a map between Banach spaces, and that \,$\Phi$\, is in fact a local diffeomorphism near \,$(u,u_y)$\,, we need to show more, however:
We need to show that the linear map \,$\Phi'((u,u_y))$\, is continuous, and has a continuous inverse. For this we consider the symplectic basis
\,$(v_k,w_k)_{k\in \Z}$\, of \,$T_{(u,u_y)}\Pot$\, mentioned above. Using its explicit representation, we can calculate its image under \,$\Phi'((u,u_y))$\,; 
by an asymptotic analysis of the components of the extended frame \,$F_\lambda$\, which comprise the \,$v_k$\, and \,$w_k$\, one can show
that the image of \,$(v_k,w_k)$\, is asymptotically close to the symplectic basis \,$(\delta \lambda_k,\delta \mu_k)_{k\in \Z}$\, of \,$T_D\Div$\,. 
The asymptotic error turns out to be sufficiently small to permit the conclusion that \,$\Phi'((u,u_y))$\, is bounded and has a bounded inverse.
Therefore \,$\Phi|\Pot_{tame}: \Pot_{tame} \to \Div_{tame}$\, is a local diffeomorphism onto an open subset of \,$\Div_{tame}$\, by the Inverse Function Theorem.

Because of the cited result on finite type divisors due to Bobenko, the image of \,$\Phi|\Pot_{tame}$\, contains all of \,$\Div_{ft,xtame}$\,, and is therefore
also dense in \,$\Div_{tame}$\,. It remains to show the injectivity of \,$\Phi|\Pot_{tame}$\,. This follows because the divisors in \,$\Div_{ft,xtame}$\, have
only one pre-image (the ``uniqueness'' part of the result by Bobenko), \,$\Div_{ft,xtame}$\, is dense in \,$\Div_{tame}$\,, and \,$\Phi|\Pot_{tame}$\, is a
local diffeomorphism. Thus the proof of Theorem~\ref{T:summary-inverse:diffeo} is concluded.

\begin{cor}
\begin{enumerate}
\item The set of potentials of finite type in \,$\Pot_{tame}$\, is dense in \,$\Pot_{tame}$\,.
\item The set of divisors \,$D\in\Div_{tame}$\, such that \,$D(x)$\, exists and is tame for all \,$x\in [0,1]$\, is open and dense in \,$\Div_{tame}$\,. 
\end{enumerate}
\end{cor}

\section{The Jacobi variety of the spectral curve}
\label{Se:summary-jacobi}

In the preceding two sections we solved the inverse problem for potentials, i.e.~we saw that a potential \,$(u,u_y)$\, can be reconstructed uniquely 
from the spectral data \,$(\Sigma,D)$\, (at least if \,$D$\, is tame). But the starting point for the present investigation was not potentials, i.e.~Cauchy data
for the sinh-Gordon equation, but rather actual simply periodic solutions \,$u$\, of the sinh-Gordon equation defined on a horizontal strip in \,$\C$\,. Therefore
we would like to understand how such an actual solution might be reconstructed from the spectral data.

Suppose \,$u: X \to \C$\, is a simply periodic solution of the sinh-Gordon equation, where \,$X\subset \C$\, is a horizontal strip in \,$\C$\,. The preceding results
have shown how to reconstruct \,$u$\, on the horizontal line through some \,$z_0\in X$\, from the spectral data \,$(\Sigma_{z_0},D_{z_0})$\, constructed from the 
monodromy at the base point \,$z_0$\, (instead of \,$z_0=0$\,, as we considered previously), or equivalently, from the spectral data (in the previous sense,
via the monodromy at the base point \,$z=0$\,) of the translated potential \,$z \mapsto u(z+z_0)$\,. One approach to the reconstruction of \,$u$\, in its entirety
is therefore via the study of how the spectral data \,$(\Sigma,D)$\, change under such a translation of the potential. 

Let us denote the monodromy of the translated potential by \,$M_{z_0}(\lambda)$\,, and the monodromy of the original potential by \,$M(\lambda)=M_{z_0=0}(\lambda)$\,. 
One can show that then \,$M_{z_0}(\lambda) = F_\lambda(z_0) \cdot M(\lambda) \cdot F_\lambda(z_0)^{-1}$\, holds (where \,$F_\lambda: X \to \SL(2,\C)$\, is the
``extended frame'', i.e.~the solution of \,$\mathrm{d}F_\lambda = \alpha_\lambda\cdot F_\lambda$\, with \,$F_\lambda(0)=\unity$\,). Therefore the eigenvalues
of the monodromy \,$M(\lambda)$\,, and hence the spectral curve \,$\Sigma$\,, does not change at all under translation.

However the eigenvectors of the monodromy, which are described by the spectral divisor \,$D$\,, of course do change under translation. It is possible to describe
the motion of the divisor points under translation by a system of differential equations for the coordinate functions \,$\lambda_k$\, and \,$\mu_k$\, of the divisor points. 
However, it turns out that this system is not locally Lipschitz continuous near infinite-type divisors when regarded on the appropriate Banach space \,$\Div$\, (locally isomorphic to 
\,$\ell^2_{-1,3} \oplus \ell^2_{0,0}$\,), so to understand the motion of the \,$\lambda_k$\, and \,$\mu_k$\, well, one needs different coordinates.

In the case of finite type divisors, it is well-known that the translations correspond to linear motions in the Jacobi coordinates of a partial desingularization of the
spectral curve (which is of finite genus in that setting). To transfer this fact to our present situation (where even the normalization of the spectral curve generally
has infinite genus), we need to construct a version of the Jacobi variety and the Abel map (hence of Jacobi coordinates) for the infinite-genus curve \,$\Sigma$\,.

Let us review the construction of the Jacobi variety for \emph{compact} Riemann surfaces 
(see for example \cite{Farkas/Kra:1992}, Section~III.6): 
Let \,$X$\, be a compact Riemann surface, say of genus \,$g \geq 1$\,, and let \,$(A_k,B_k)_{k=1,\dotsc,g}$\,
be a canonical homology basis of \,$X$\,, i.e.~\,$(A_k,B_k)$\, is a basis of the homology group \,$H_1(X,\Z)$\, 
with the intersection properties \,$A_k\times B_\ell = \delta_{k\ell}$\, (Kronecker delta), \,$A_k\times A_\ell = 0 = B_k \times B_\ell$\,
for \,$k,\ell\in \{1,\dotsc,g\}$\,. 
Then there exists a canonical basis \,$(\omega_k)_{k=1,\dotsc,g}$\, of the vector space \,$\Omega(X)$\, of holomorphic 1-forms on \,$X$\, 
that is dual to \,$(A_k)$\, in the sense that \,$\int_{A_k} \omega_\ell = \delta_{k\ell}$\, holds. To any given positive divisor \,$D=\{P_1,\dotsc,P_g\}$\, of degree \,$g$\, on \,$X$\,, we then associate
the quantity
$$ \wt{\vi}(D) := \left( \sum_{k=1}^g \int_{P_0}^{P_k} \omega_\ell \right)_{\ell=1,\dotsc,g} \;\in\; \C^g\;, $$
where \,$P_0\in X$\, is the ``origin point'', which we hold fixed. 
Because these integrals depend on the homology class of the paths of integration from \,$P_0$\, to \,$P_k$\, we choose, the quantity \,$\wt{\vi}(D)$\, is only defined modulo the \emph{period lattice}
$$ \Gamma := \left\langle \left( \int_{A_k} \omega_\ell \right)_{\ell=1,\dotsc,g} \;,\; \left( \int_{B_k} \omega_\ell \right)_{\ell=1,\dotsc,g} \right\rangle_{\Z} \subset \C^g \; . $$
Thus we obtain the \emph{Jacobi variety} \,$\Jac(X) := \C^g / \Gamma$\, of \,$X$\, and by projecting the values of \,$\wt{\vi}$\, onto \,$\Jac(X)$\, the \emph{Abel map} 
\,$\vi: \mathrm{Div}_g(X) \to \Jac(X)$\,, 
where \,$\mathrm{Div}_g(X)$\, denotes the space of positive divisors of degree \,$g$\, on \,$X$\,. 

We need to generalize this construction for the spectral curve \,$\Sigma$\,. In particular we need to deal with the fact that \,$\Sigma$\, is not compact,
and its homology group is generally infinite dimensional. 
In particular the space \,$\C^g$\, occurring in the treatment of the compact case as universal cover of the Jacobi variety will need to be replaced by
an infinite-dimensional Banach space, and the sum defining the Abel map will be an infinite sum. To ensure its convergence, we will need to impose
a condition on the divisors we admit for the Abel map, and this condition is precisely the asymptotic condition for the space \,$\Div$\, given in 
Definition~\ref{D:summary-remon:asymptotic-nonspecial}(1).

For the purposes of this summary paper, we will ignore the second complication that would need to be addressed, namely that \,$\Sigma$\,
can have singularities. In other words, in the present
and the following section we will always suppose that \,$\Sigma$\, does not have any singularities and thus is a Riemann surface.

Let \,$\Sigma$\, be the spectral curve corresponding to some potential \,$(u,u_y)\in \Pot$\, satisfying the above hypothesis. We enumerate the branching points
of \,$\Sigma$\, by two sequences \,$(\vkap_{k,1})_{k\in \Z}$\, and \,$(\vkap_{k,2})_{k\in \Z}$\, as in Theorem~\ref{T:summary-asymp:spectral}(1). Because \,$\Sigma$\,
is a Riemann surface, we have \,$\vkap_{k,1}\neq \vkap_{k,2}$\, for all \,$k\in \Z$\,. We now fix a homology basis for \,$\Sigma$\,: 
For \,$k\in \Z$\, we let \,$A_k$\, be a small, non-trivial cycle in \,$\Sigma$\,  that 
encircles the pair of branch points \,$\vkap_{k,1}$\, and \,$\vkap_{k,2}$\, of \,$\Sigma$\,, but no other branch points. For \,$k\neq 0$\, there exists another
cycle \,$B_k$\, that encircles the branch points \,$\vkap_{k,1}$\, and \,$\vkap_{k=0,1}$\,, and no others. 
A final cycle \,$B_0$\, comes from the observation
that \,$\sqrt{\lambda}$\, is a global parameter on \,$\Sigma$\, away from the branch points. Because
the Riemann surface associated to \,$\sqrt{\lambda}$\, has branch points in \,$\lambda=0$\, and \,$\lambda=\infty$\,, we see that \,$\Sigma$\, also has branch points there,
and thus there is another non-trivial cycle \,$B_0$\, that encircles these two branch points. We choose the orientation of all these cycles so that their intersection
numbers satisfy
\,$A_k\times A_\ell = 0 = B_k\times B_\ell$\, and \,$A_k\times B_\ell = \delta_{k\ell}$\, (Kronecker delta) for all \,$k,\ell$\,. 
Then \,$(A_k,B_k)_{k\in\Z}$\, is a canonical basis of the homology of \,$\Sigma$\,.

One can show that \,$\Sigma$\, is parabolic in the sense of \textsc{Ahlfors} and \textsc{Nevanlinna} (see for example \cite[Chapter 1]{Feldman/Knoerrer/Trubowitz:2003}),
and from this fact it follows that there exists a basis \,$(\omega_n)_{n\in \Z}$\, of the space \,$\Omega(\Sigma) \cap L^2(\Sigma,T^*\Sigma)$\, of square-integrable,
holomorphic 1-forms on \,$\Sigma$\, that is dual to the canonical basis of the homology \,$(A_k,B_k)$\, in the sense that
\,$\int_{A_k} \omega_n = \delta_{k,n}$\, (Kronecker delta) holds. However, as it is described in \cite[Chapter~17 and the first half of Chapter~18]{Klein:2015},
for our specific situation with \,$\Sigma$\, being a spectral curve, it is possible to give an explicit description of the \,$\omega_n$\, as a linear combination
of infinite products. This explicit description is useful because by its investigation one can show that the \,$\omega_n$\, show a steeper descent towards zero
for \,$\lambda\to 0$\, and \,$\lambda\to\infty$\, than is expressed by the fact that they are square-integrable alone. This steeper asymptotic behavior turns 
out to be crucial for the construction of the Abel map. 

We consider the periods corresponding to the \,$\omega_n$\,, i.e.~for \,$k,n\in \Z$\, we let
$$ \alpha^{[k]}_n := \int_{A_k} \omega_n = \delta_{k,n} \qmq{and} \beta^{[k]}_n := \int_{B_k} \omega_n \; . $$
Using the explicit description of the \,$\omega_n$\, that was mentioned above, it can be shown (see \cite[Theorem~18.7(1)]{Klein:2015}) that 
the \,$\beta^{[k]}_n$\, satisfy the asymptotic property
$$ \bigr(\beta^{[k]}_n \cdot (\vkap_{n,1}-\vkap_{n,2})\bigr)_{n\in \Z} \in \ell^2_{-1,3} \qmq{for every \,$k\in \Z$\,.} $$
Thus we are led to consider the Banach space
$$ \wt{\Jac}(\Sigma) := \Mengegr{(a_n)_{n\in \Z}}{a_n \cdot (\vkap_{n,1}-\vkap_{n,2}) \in \ell^2_{-1,3}} $$
with the norm 
$$ \|a_n\|_{\wt{\Jac}(\Sigma)} := \left\|a_n \cdot (\vkap_{n,1}-\vkap_{n,2})\right\|_{\ell^2_{-1,3}} \qmq{for \,$(a_n) \in \wt{\Jac}(\Sigma)$\,.} $$
Then we have \,$(\alpha^{[k]}_n)_{n\in\Z}, (\beta^{[k]}_n)_{n\in \Z} \in \wt{\Jac}(\Sigma)$\, for every \,$k\in \Z$\,. For this reason we use \,$\wt{\Jac}(\Sigma)$\,
in the place of \,$\C^g$\, in the construction of the Jacobi variety. It should be mentioned that \,$\wt{\Jac}(\Sigma)$\, is a Banach space only under our
hypothesis that \,$\Sigma$\, is regular. In the more general case where \,$\Sigma$\, has singularities, \,$\|\,\cdot\,\|_{\wt{\Jac}(\Sigma)}$\, is only a semi-norm,
and thus \,$\wt{\Jac}(\Sigma)$\, becomes only a topological vector space with the induced (non-Hausdorff) topology. 

We now fix an asymptotic divisor \,$D^o$\, on \,$\Sigma$\, which will serve as the origin divisor for the construction of the Abel map. 
To ensure that
the infinite sum that will define the Abel map converges, we need to restrict the integration paths we consider. 
For this purpose we denote by \,$\mathfrak{C}$\, the set of sequences \,$(\gamma_k)_{k\in \Z}$\, where each \,$\gamma_k$\, is a curve in \,$\Sigma$\, running from
a point \,$(\lambda_k^o,\mu_k^o) \in \Sigma$\, to another point \,$(\lambda_k,\mu_k) \in \Sigma$\,, such that \,$(\lambda_k^o,\mu_k^o)_{k\in \Z}$\, equals the support of \,$D^o$\, and 
the divisor \,$D$\, with support \,$(\lambda_k,\mu_k)_{k\in \Z}$\, is another asymptotic divisor; 
moreover for large \,$|k|$\, the curve \,$\gamma_k$\, winds around no branch points of \,$\Sigma$\, but \,$\vkap_{k,1}$\, and \,$\vkap_{k,2}$\,,
and there is a number \,$m_\gamma\in \N$\, (depending on \,$\gamma$\, but not on \,$k$\,) so that the winding number of any \,$\gamma_k$\, around any branch point of \,$\Sigma'$\,
is at most \,$m_\gamma$\,. In this situation we call \,$D$\, the divisor induced by the sequence of curves \,$(\gamma_k)_{k\in \Z}$\,. 
Every asymptotic divisor on \,$\Sigma$\, is induced by some \,$(\gamma_k) \in \mathfrak{C}$\, in this sense.

We let \,$\Gamma$\, be the abelian group corresponding to the periods of all closed loops in \,$\mathfrak{C}$\,, i.e.~
$$ \Gamma := \left. \left\{ \sum_{k\in \Z} (a_k\,\alpha^{[k]}+b_k\,\beta^{[k]}) \right| \begin{matrix} a_k,b_k\in \Z \\ \exists N,m\in \N\;\forall k\in \Z,|k|>N\;:\; |a_k|\leq m,\;b_k=0 \end{matrix} \right\} \; . $$
Then \,$\Gamma$\, is an abelian subgroup of \,$\wt{\Jac}(\Sigma)$\,. However, \,$\Gamma$\, is not a discrete subset of \,$\wt{\Jac}(\Sigma)$\, (\,$0$\, is an
accumulation point of \,$\Gamma$\,). Despite this fact we call \,$\Gamma$\, the \emph{period lattice} of \,$\Sigma$\,, and
we call the topological quotient space \,$\Jac(\Sigma) := \wt{\Jac}(\Sigma)/\Gamma$\, the \emph{Jacobi variety} of \,$\Sigma$\,. We denote the
canonical projection map by \,$\pi: \wt{\Jac}(\Sigma) \to \Jac(\Sigma)$\,.

The situation with \,$\Gamma$\, being non-discrete is similar to the one encountered by \hbox{\textsc{McKean}} and \textsc{Trubowitz} in \cite{McKean/Trubowitz:1976} concerning the Jacobi
variety for the integrable system associated to Hill's operator: 
There the period lattice is also not discrete in the respective Banach space, and the Jacobi variety is compact (topologically, it is a product of infinitely many circles)
and therefore does not carry the structure of an infinite dimensional manifold, see the discussion in \cite{McKean/Trubowitz:1976}, p.~154. 

The following statement is shown via a detailed asymptotic analysis of the integral \,$\int_{\gamma_k}\omega_n$\, for \,$(\gamma_k) \in \mathfrak{C}$\,,
again involving the result on the asymptotic descent of the \,$\omega_n$\, near \,$\lambda=0$\, and \,$\lambda=\infty$\,,
see \cite[Sections 16, 17, and Theorem~18.5]{Klein:2015}:
For \,$(\gamma_k)_{k\in \Z} \in \mathfrak{C}$\, and \,$n\in \Z$\,, the infinite sum \,$\sum_{k\in \Z} \int_{\gamma_k} \omega_n$\, converges absolutely in \,$\C$\,,
and if we define
$$ \wt{\vi}_n : \mathfrak{C}\to \C,\; (\gamma_k)_{k\in \Z} \mapsto \sum_{k\in \Z} \int_{\gamma_k} \omega_n \;, $$
we have
$$ \bigr(\; \wt{\vi}_n((\gamma_k)) \;\bigr)_{n\in \Z} \in \wt{\Jac}(\Sigma) \; . $$
We thus define \,$\wt{\vi}: \mathfrak{C} \to \wt{\Jac}(\Sigma),\; (\gamma_k)_{k\in \Z} \mapsto \bigr(\; \wt{\vi}_n((\gamma_k)) \;\bigr)_{n\in \Z}$\,.

We denote the space of asymptotic divisors on \,$\Sigma$\, by \,$\Div(\Sigma)$\, and we let \,$\tau: \mathfrak{C} \to \Div(\Sigma)$\, be the 
surjective map that associates to each \,$(\gamma_k)_{k\in \Z} \in \mathfrak{C}$\, the divisor induced by \,$(\gamma_k)$\,. 
Then there exists one and only one map \,$\vi: \Div(\Sigma) \to \Jac(\Sigma)$\, with \,$\vi\circ \tau = \pi \circ \wt{\vi}$\,,
i.e.~so that the following diagram commutes (see \cite[Theorem~18.7(2)]{Klein:2015}):
\begin{equation*}
\begin{minipage}{5cm}
\begin{xy}
\xymatrix{
\mathfrak{C}_{D^o} \ar[r]^{\wt{\vi}} \ar[d]_{\tau} & \wt{\Jac}(\Sigma) \ar[d]^{\pi} \\
\Div(\Sigma) \ar[r]_{\vi} & \Jac(\Sigma)
}
\end{xy}
\end{minipage}
\end{equation*}
We call \,$\vi$\, the \emph{Abel map} of \,$\Sigma$\,. It is clear that change of the origin divisor \,$D^o$\, corresponds to a linear transformation
of \,$\vi$\, (see \cite[Theorem~18.7(3)]{Klein:2015}).

The space \,$\wt{\Jac}(\Sigma)$\, plays the role of a tangent space for the Jacobi variety \,$\Jac(\Sigma)$\,. In our setting where the period lattice \,$\Gamma$\, is not discrete,
the tangent space of \,$\Jac(\Sigma)$\, is not unique however, and similarly as it is the case for Hill's equation as studied by \textsc{McKean} and \textsc{Trubowitz} in \cite{McKean/Trubowitz:1976}, 
we need to pass to a larger tangent space so that the flow of translations of the potential (which we will study via the Jacobi variety in the following section)
is tangential to \,$\Jac(\Sigma)$\,. This corresponds to a larger space of curve tuples \,$\mathfrak{C}^{(1)}$\,
and a larger Banach space \,$\wt{\Jac}^{(1)}(\Sigma)$\,.
In fact McKean and Trubowitz construct in \cite{McKean/Trubowitz:1976} an entire ascending family of tangent spaces for their Jacobi variety, which correspond to the higher flows
of the integrable system associated with Hill's equation. In our setting we cannot define more than the first extension described in the following proposition, because
our potentials are only once weakly differentiable, in contrast to the infinitely differentiable potentials in \cite{McKean/Trubowitz:1976}.

Explicitly,
let \,$\mathfrak{C}^{(1)}$\, be the set of sequences
\,$(\gamma_k)_{k\in \Z}$\, where each \,$\gamma_k$\, is a curve in \,$\Sigma$\, running from
a point \,$(\lambda_k^o,\mu_k^o)\in \Sigma$\, to another point \,$(\lambda_k,\mu_k) \in \Sigma$\,, such that \,$(\lambda_k^o,\mu_k^o)_{k\in \Z}$\, equals the support of \,$D^o$\, and
the divisor \,$D$\, with support \,$(\lambda_k,\mu_k)_{k\in \Z}$\, is another asymptotic divisor; moreover for large \,$|k|$\, the curve \,$\gamma_k$\, winds around
no branch points of \,$\Sigma$\, but \,$\vkap_{k,1}$\, and \,$\vkap_{k,2}$\,  
and there is a number \,$m_\gamma\in \N$\, (depending on \,$\gamma$\, but not on \,$k$\,) so that the winding number of any \,$\gamma_k$\, around any branch point or puncture of \,$\Sigma'$\,
is at most \,$m_\gamma\cdot |k|$\,. Then we define the Banach space
$$ \wt{\Jac}^{(1)}(\Sigma) := \Mengegr{(a_n)_{n\in \Z}}{a_n \cdot (\vkap_{n,1}-\vkap_{n,2}) \in \ell^2_{-2,2}} $$
and the maps
$$ \wt{\vi}_{n}^{(1)} : \mathfrak{C}^{(1)} \to \C,\;(\gamma_k)_{k\in \Z} \mapsto \sum_{k\in \Z} \int_{\gamma_k} \omega_n \;. $$
Essentially in the same way as above one shows (see \cite[Proposition~18.10]{Klein:2015}) that the sum defining \,$\wt{\vi}_n^{(1)}$\, is still absolutely convergent and that 
\,$\bigr(\wt{\vi}_n^{(1)}((\gamma_k))\bigr)_{n\in \Z} \in \wt{\Jac}^{(1)}(\Sigma)$\, holds for any \,$(\gamma_k)_{k\in \Z} \in \mathfrak{C}^{(1)}$\,. 
Thus we obtain an extended Jacobi coordinate map
$$ \wt{\vi}^{(1)}: \mathfrak{C}^{(1)} \to \wt{\Jac}^{(1)}(\Sigma),\; (\gamma_k)_{k\in \Z} \mapsto \bigr(\wt{\vi}_n^{(1)}((\gamma_k))\bigr)_{n\in \Z} \; . $$
Clearly \,$\wt{\Jac}(\Sigma) \subset \wt{\Jac}^{(1)}(\Sigma)$\, and \,$\wt{\vi}^{(1)}|\mathfrak{C} = \wt{\vi}$\, holds.

\section{Translations of divisors, and the asymptotic behavior of spectral data for simply periodic solutions}
\label{Se:summary-xytrans}

We are now ready to describe the motion of the points 
of asymptotic divisors under translation (in the sense explained at the beginning of Section~\ref{Se:summary-jacobi}) in terms of Jacobi coordinates.
For this purpose we continue to use the notations of the previous setting. For spectral data \,$(\Sigma,D)$\, of a simply periodic solution \,$u$\, of the sinh-Gordon equation,
we denote by \,$D(x)$\, resp.~\,$D(y)$\, the spectral divisor of the solution \,$u$\, translated in \,$x$-direction resp.~in \,$y$-direction
(also see the discussion at the beginning of the previous section). 
For the construction of the Abel map on the spectral curve \,$\Sigma$\,, we fix the origin divisor as \,$D^o := D(0)=D$\,. 

In the sequel we will look at the derivatives \,$\tfrac{\partial \vi_n}{\partial x}$\, and \,$\tfrac{\partial \vi_n}{\partial y}$\, of the \,$n$-th Jacobi coordinate \,$\vi_n$\,. 
For these derivatives to make sense,
we need to define Jacobi coordinates \,$\vi_n$\, of \,$D(x)$\, resp.~\,$D(y)$\, at least for small \,$|x|$\, resp.~\,$|y|$\, for all \,$n\in \Z$\,. For this purpose
we write the support of \,$D(x)$\, as \,$(\lambda_k(x),\mu_k(x))_{k\in \Z}$\, for small \,$|x|$\, and then consider for fixed \,$x$\, and all \,$k\in \Z$\, the curve 
\,$\gamma_{x,k}: [0,x] \to \Sigma,\; t\mapsto (\lambda_k(t),\mu_k(t))$\,. Because the spectral map \,$\Pot\to\Div$\, is asymptotically close to the Fourier transform of the potential, \,$\gamma_{x=1,k}$\, winds \,$|k|$\, times around
the pair of branch points \,$\vkap_{k,1}$\,, \,$\vkap_{k,2}$\, for \,$|k|$\, large; it follows that we do not have \,$(\gamma_{x,k}) \in \mathfrak{C}$\,,
but we do have \,$(\gamma_{x,k}) \in \mathfrak{C}^{(1)}$\,. Therefore we can define Jacobi coordinates for the translation in \,$x$-direction in the vicinity of \,$D(0)$\, by
$$ \vi_n(x) := \wt{\vi}_n^{(1)}(\gamma_{x,k}) \qmq{for \,$n\in \Z$\,.} $$
A similar construction applies for the translation in \,$y$-direction; here it is relevant that for large \,$|k|$\, the divisor
point \,$(\lambda_k(y),\mu_k(y))$\, remains close to the pair of branch points \,$\vkap_{k,1}$\,, \,$\vkap_{k,2}$\, 
for sufficiently small \,$|y|$\, because the asymptotic estimates then apply to the translated potentials uniformly.

We denote by \,$\tfrac{\partial \vi_n}{\partial x}$\,
resp.~\,$\tfrac{\partial \vi_n}{\partial y}$\, the derivative of the Jacobi coordinate \,$\vi_n(x)$\, resp.~\,$\vi_n(y)$\, with respect to \,$x$\, resp.~\,$y$\,. 

The following theorem expresses that like in the finite-type setting, also in our present situation where the spectral curve is of infinite geometric
genus, the translations of the divisor correspond to linear motions in the Jacobi variety. 

\begin{thm}
\label{T:summary-xytrans:jacobi}
There exist sequences \,$a_n^x,a_n^y \in \ell^2_{-1,-1}$\, (dependent only on the spectral curve \,$\Sigma$\,) so that 
under translation of the potential \,$u$\, in the direction of \,$x$\, resp.~\,$y$\,, the Jacobi coordinates \,$\vi_n$\, (\,$n\in \Z$\,) follow the differential equations
\begin{align*}
\frac{\partial \vi_n}{\partial x} & = n + a_n^x \;, \\
\frac{\partial \vi_n}{\partial y} & = -i|n| + a_n^y \; . 
\end{align*}
Moreover we have \,$a_n^x=0$\, for \,$|n|$\, large, and for every \,$n\in \Z$\,, \,$\tfrac{\partial \vi_n}{\partial x}$\, corresponds
to a member of the period lattice, i.e.~there exists a cycle \,$Z_n$\, of \,$\Sigma$\, so that \,$\tfrac{\partial \vi_n}{\partial x}
= \int_{Z_n} \omega_n$\, holds.
\end{thm}

The statement that \,$\tfrac{\partial \vi_n}{\partial x}$\, is a member of the period lattice of \,$\Sigma$\, corresponds to the fact that the solution \,$u$\,
of the sinh-Gordon equation is periodic in the \,$x$-direction. This is why there is no analogous statement for \,$\tfrac{\partial \vi_n}{\partial y}$\, in general.

At the heart of the proof of Theorem~\ref{T:summary-xytrans:jacobi}, which is detailed in \cite[Section~19]{Klein:2015},
is a general construction of linear flows in the Picard variety of a Riemann surface \,$X$\,
(the space of isomorphy classes of line bundles on \,$X$\,) known as the \emph{Krichever construction}. In fact it turns out that the vector fields
\,$\tfrac{\partial \vi_n}{\partial x}$\, and \,$\tfrac{\partial \vi_n}{\partial y}$\, can be constructed on \,$\Sigma$\, via the Krichever construction
by marking the points \,$\lambda=0$\, and \,$\lambda=\infty$\, and prescribing suitable Laurent series with poles of order 1 around these points.

\bigskip

Because solutions of the sinh-Gordon equation are real analytic on the interior of their domain of definition, we expect that spectral data \,$(\Sigma,D)$\, 
corresponding to simply periodic solution of the sinh-Gordon equation on a horizontal strip of positive height to have a far better asymptotic behavior
than the relatively mild asymptotic law for spectral data of potentials \,$(u,u_y)\in \Pot$\, with merely \,$u \in W^{1,2}([0,1])$\, and \,$u_y \in L^2([0,1])$\, 
that was found in Theorem~\ref{T:summary-asymp:spectral}. 
More specifically, we expect both the distance of branch points \,$\vkap_{k,1}-\vkap_{k,2}$\, of the spectral curve \,$\Sigma$\, and the distance of the corresponding spectral divisor points
to the branch points to fall off exponentially for \,$k\to\pm \infty$\,. The following theorem shows that our expectations are correct:

\begin{thm}
Let \,$y_0,\eps>0$\,, \,$X = \Mengegr{z\in \C}{|\IM(z)|< y_0+\eps}$\, the horizontal strip in \,$\C$\, of height \,$2(y_0+\eps)$\,, and \,$u:X \to \C$\, be a 
simply periodic solution of the
sinh-Gordon equation \,$\Delta u + \sinh(u)=0$\,.
We let \,$\Sigma$\, be the spectral curve corresponding to \,$u$\, (with branch points \,$\vkap_{n,\nu}$\,, and \,$\vkap_{n,*}:=\tfrac12(\vkap_{n,1}+\vkap_{n,2})$\,) 
and let \,$D := \{(\lambda_n,\mu_n)\}_{n\in \Z}$\, be the spectral divisor of \,$u$\, with the starting point \,$z_0=0$\,. 

Then there exists a constant \,$C>0$\, and a sequence \,$(s_n)_{n\in \Z} \in \ell^2_{0,0}$\, of real numbers so that 
\begin{align*} 
|\vkap_{n,1}-\vkap_{n,2}| & \leq C\,e^{-2\pi\,(1-s_n)\,|n|\,y_0}  \;,\\
|\lambda_n-\vkap_{n,*}| & \leq C\,e^{-2\pi\,(1-s_n)\,|n|\,y_0} \;,  \\
\qmq{and} |\mu_n - (-1)^n| & \leq C\,e^{-\pi\,(1-s_n)\,|n|\,y_0} \;. 
\end{align*}
\end{thm}

The proof of this theorem is described in \cite[Section~20]{Klein:2015}. It is based on the description of the flow of the Jacobi coordinates under translations
in Theorem~\ref{T:summary-xytrans:jacobi}, in conjunction with a careful analysis of the asymptotic behavior of the Abel map.

\bigskip
\end{document}